\documentclass[12pt]{article}
\usepackage[final]{epsfig}
\usepackage{graphics}
\usepackage{amsmath}
\usepackage{amsfonts}
\usepackage{latexsym}
\usepackage{amssymb}
\usepackage{graphicx}
\usepackage{url}
\usepackage{epstopdf}
\usepackage{hyperref}
\usepackage{color}
\usepackage{marginnote}
\usepackage{a4wide}

\newtheorem{lemma}{Lemma}[section]

\newtheorem{remark}[lemma]{Remark}
\newtheorem{example}[lemma]{Example}
\newtheorem{theorem}{Theorem}

\newtheorem{corollary}[lemma]{Corollary}

\newcommand{\g}{{\gamma}}
\newcommand{\eps}{{\varepsilon}}
\newcommand{\proofend}{$\Box$\bigskip}

\newcommand{\R}{{\mathbb R}}

\def\proof{\paragraph{Proof.}}

\begin{document}

\title{Outer length billiards on a large scale}

\author{Peter Albers\footnote{
Mathematisches Institut,
Universit\"at Heidelberg,
69120 Heidelberg,
Germany;
peter.albers@uni-heidelberg.de}
 \and 
 Lael Edwards-Costa\footnote{
Department of Mathematics,
Pennsylvania State University,
University Park, PA 16802,
USA;
lael.costa@psu.edu}
 \and
 Serge Tabachnikov\footnote{
Department of Mathematics,
Pennsylvania State University,
University Park, PA 16802,
USA;
tabachni@math.psu.edu}
} 

\date{\today}

\maketitle


\section{Introduction} \label{sect:three}

The billiard dynamical system describes the motion of a point mass confined in a domain and experiencing elastic collision with its boundary. 
We consider the case when the billiard table is bounded by a smooth closed quadratically convex planar curve (an oval).

\begin{figure}[ht]
\centering
\includegraphics[width=1.8in]{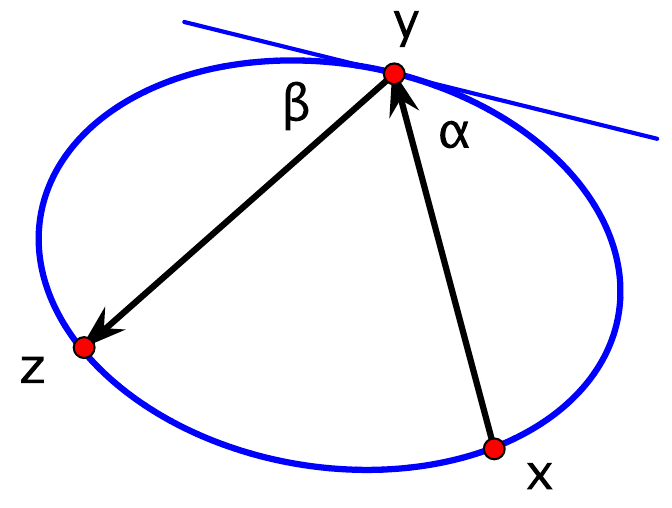}\quad
\includegraphics[width=1.8in]{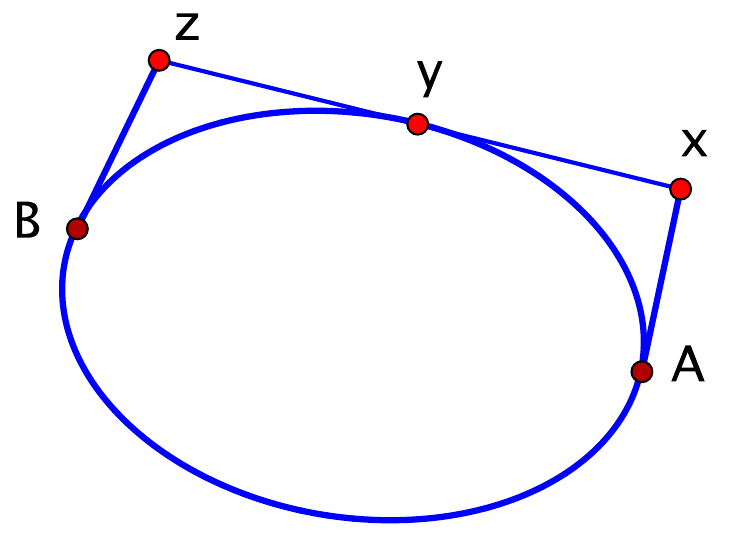}\quad
\includegraphics[width=1.8in]{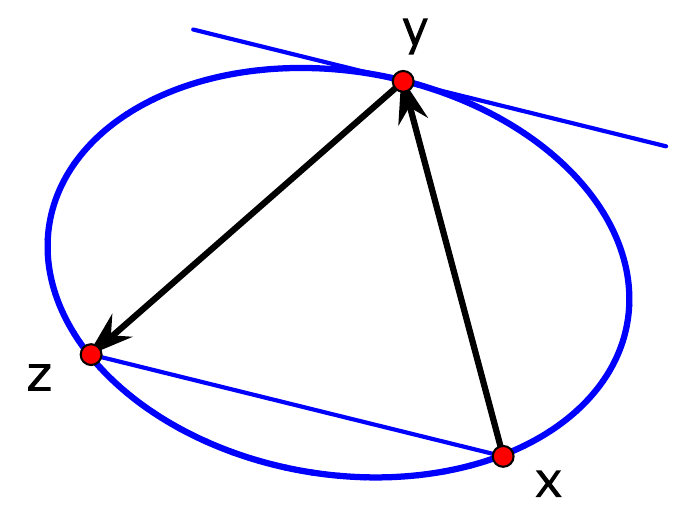}
\caption{Three billiards.}	
\label{three}
\end{figure}

The billiard reflection law, the angle of incidence equals the angle of reflection, has a variational origin. Consider the left Figure \ref{three}. Fix points $x$ and $z$ and let $y$ vary on the curve $\g$. Then $\alpha=\beta$ if and only if the length $|xy|+|yz|$ is extremal as a function of $y\in \g$, see, e.g., \cite{Ta}. 
Thus one may describe this system as the ``inner length billiard".   

Now consider the middle Figure \ref{three}. It depicts the outer billiard system, an area preserving map of the exterior of an oriented oval $\g$: the map send $x$ to $z$ if the segment is tangent to the with the same orientation, and the tangency point $y$ is the midpoint of $xz$. 

This reflection law also has a variational description. Fix points $A$ and $B$, and let the line $xz$ be a variable tangent line to $\g$ that intersects the tangent lines at point $A$ and $B$ at points $x$ and $z$ respectively. Then $|xy|=|yz|$ if and only if the area bounded by the segments $Ax, xz,zB$ and the arc $AyB$ is extremal as a function of $y$. See, e.g., \cite{DT} for a survey of outer billiards. Outer billiards may be described as the ``outer area billiards".

Now consider the right Figure \ref{three}. This is symplectic billiard, introduced in \cite{AT}, also see \cite{BBN1,BB,BBN2}. Like the conventional billiard, this is a map of oriented chords of an oval $\g$. The chord $xy$ is mapped to $yz$ when the tangent line to $\g$ at $y$ is parallel to $xz$.

This non-local reflection law is also variational: if points $x$ and $z$ are fixed and point $y\in \g$ varies, then $xz \parallel T_y \g$ if and only if the area of the triangle $xyz$ is extremal as  a function of $y$. Thus symplectic billiard is the ``inner area billiard".\footnote{Like the conventional billiard, the outer and symplectic billiards can be defined in multi-dimensional setting: the ambient space is linear symplectic.}

One cannot help wondering about the missing fourth,  ``outer length", billiard. This ``4th billiard" was recently  introduced in \cite{AT1} and already studied in \cite{BBF} and \cite{BBS}. See also the survey \cite{BM}. 
We continue the study of this billiard-like system in the present paper. 

\begin{figure}[ht]
\centering
\includegraphics[width=.45\textwidth]{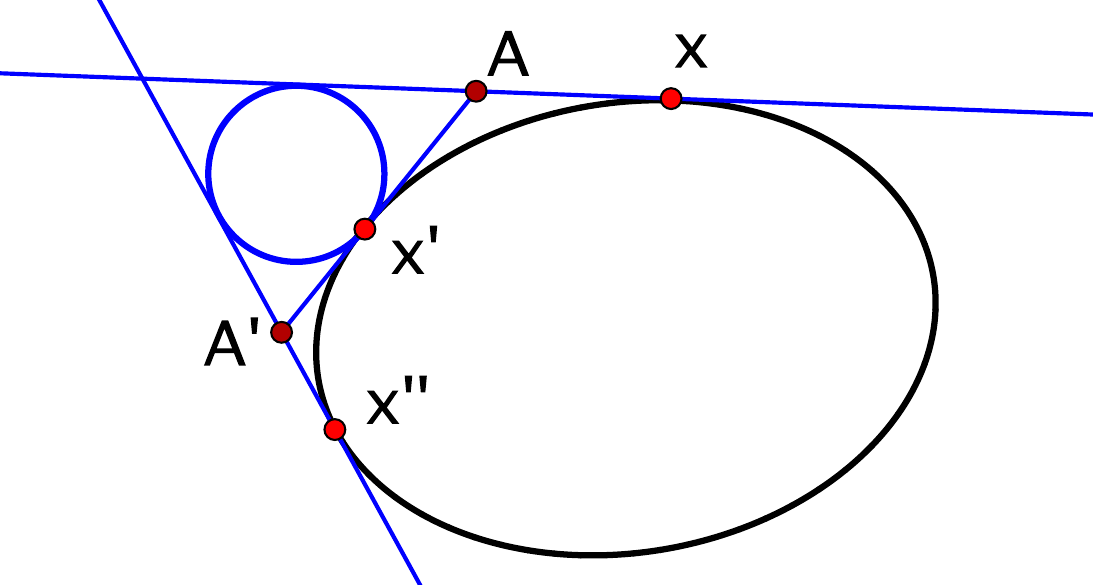}
\caption{The outer length billiard.}	
\label{fourt}
\end{figure}

The geometric definition of the outer length billiard reflection $F$ is rather involved, see Figure \ref{fourt}. For a deduction of this geometric construction from the variational definition (length extremizing), see \cite{AT1}.

Let $A$ be a point outside of 
$\g$, and let $Ax'$ and $Ax$ be the positive and negative tangent segments to $\g$ (the sign given by the orientation of the oval). Consider the circle tangent to $\g$ at point $x'$, tangent to the line $Ax$, and lying on the right of the ray $Ax'$. Then $F(A)=A'$ is defined as the intersection point of the line $Ax'$ and the common tangent line of the circle and $\g$.

Note that this construction still makes  sense if $\g$ is a convex polygon -- as long as the point $A$ does not lie on the extension of a side. This assumption is satisfied in a set of full measure, thus 
one can study polygonal outer length billiards as well. See the accompanying paper \cite{EC} for a study of the polygonal case. 

\begin{figure}[ht]
\centering
\includegraphics[width=.45\textwidth]{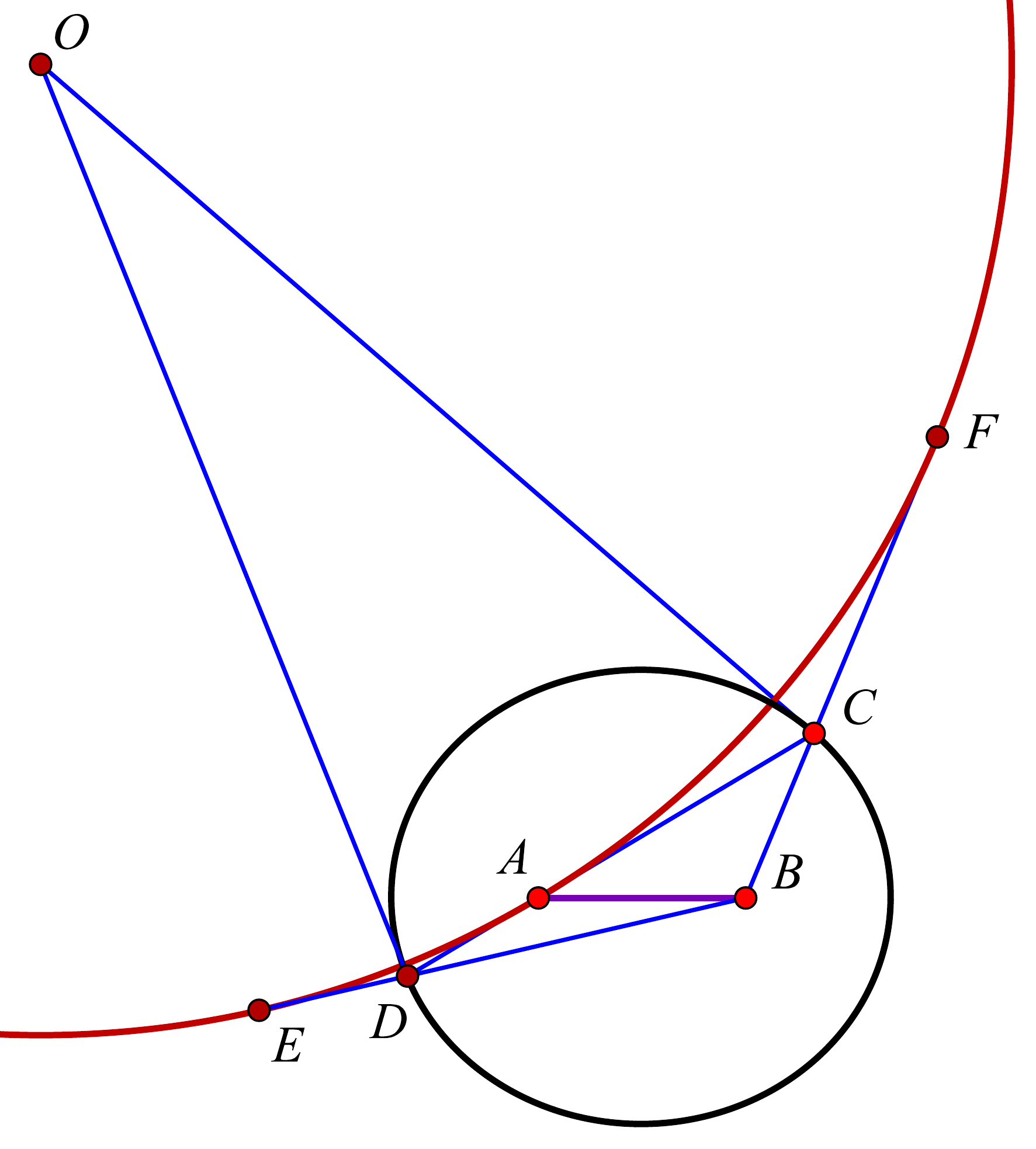}
\caption{Outer length billiard map on a segment $AB$: one has $F(C)=D$.}	
\label{segment}
\end{figure}

For example, Figure \ref{segment} depicts the outer length billiard on a segment. A confocal ellipse is an invariant curve. Indeed, since two tangents to a circle from a point are equal, one has
$$
BC+CA=BC+CF=BD+DE=BD+DA,
$$ 
hence points $C$ and $D$ lie on an ellipse with foci $A$ and $B$.\footnote{In contrast, the second iteration of the outer (area) billiard map on a segment is parallel translation in the direction of the segment through distance equal to twice its length.} See Section \ref{sect:ell} for the outer length billiard maps on ellipses.

The goal of this paper is two-fold. Section \ref{sect:def} provides some foundational results, presented  without proofs in \cite{AT1}: the generating function and the invariant area form of the outer length billiard map. We also provide various comment on the important fact: the 
the outer length billiard map on ellipses is completely integrable, confocal ellipses being its invariant curves. This fact was proved in \cite{BBF}. 

The second goal of the paper is motivated by our recent study of the dynamical behavior of the outer symplectic billiard map ``at infinity", that is, very far away from the billiard table, in \cite{ACT} (see \cite{Ta96} for a precursor of this result, in 2D). The main result of \cite{ACT} is that the second iteration of the map is well approximated by the flow of a Hamiltonian vector field whose homogeneous of degree 1 Hamiltonian is determined by the shape of the billiard table. 

\begin{figure}[ht]
\centering
\includegraphics[width=.4\textwidth]{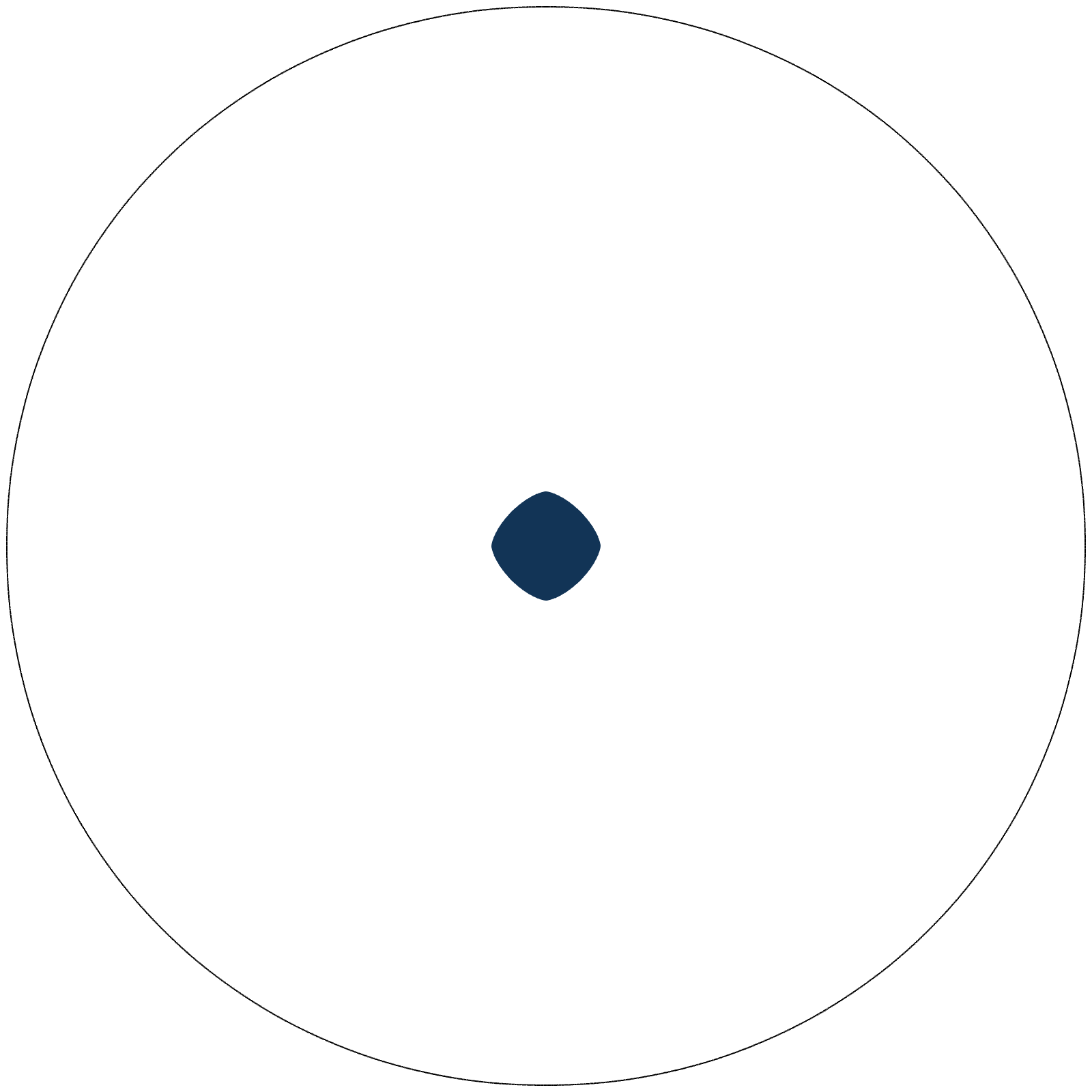}\qquad\qquad
\includegraphics[width=.4\textwidth]{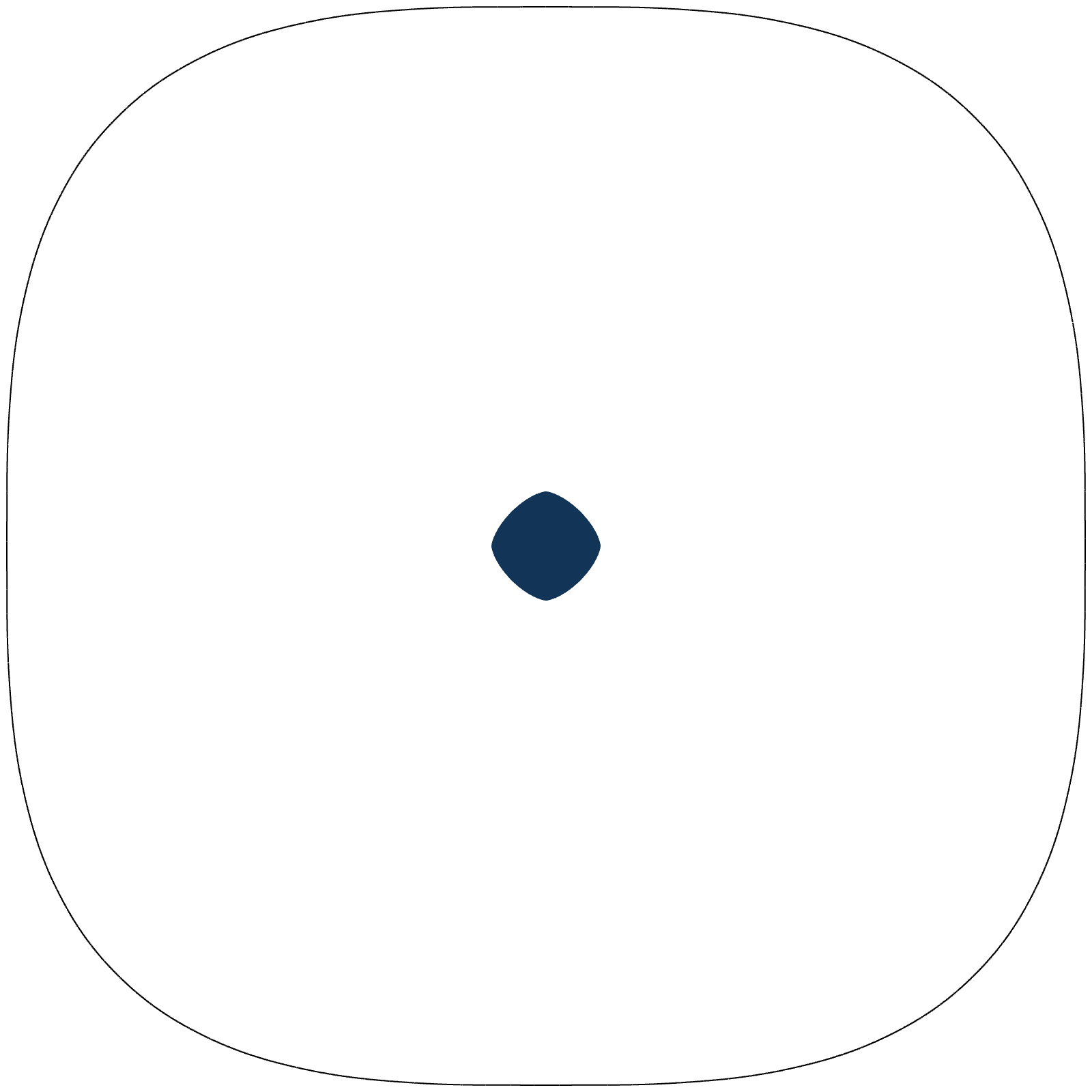}
\caption{Left: the orbit of a point ``at infinity": the table is a circle in the $L_p$-norm with $p=3/2$.
Right: the orbit of the centers of the auxiliary circles ``at infinity": the orbit a circle in the $L_p$-norm with $p=3$, the polar dual to the table.}	
\label{triinfty}
\end{figure}

Take a point very far away from the billiard table, the oval $\g$, and consider its orbit under the second iteration of the outer length billiard map. To fit into a computer screen,  rescale so that the the billiard table appears as almost a point. Then the orbit appears as a circle, see Figure \ref{triinfty}, and the discrete time motion appears to be continuous. On can experimentally evaluate the speed of this motion at a point of this large circle parameterized as  $r(\cos \alpha,\sin\alpha)$ (with an origin chosen inside $\g$): this speed is proportional to $w_\gamma(\alpha)$, the width of the oval $\g$ in the direction $\alpha$. 
This phenomenon is studied and explained in Section \ref{sect:dyninfty}. See Theorem \ref{thm:main1} for a precise formulation. 

We also show that if $\g$ is sufficiently smooth then the outer length billiard map possesses invariant curves arbitrarily far away from the table $\g$. This is an analog of a similar observation by J. Moser concerning outer (area) billiards. As a consequence, one has stability: the orbits do not escape to infinity.

Next we consider the orbits at infinity of the centers of the auxiliary circles involved in the definition of the outer length billiard map. To see these orbits one needs to further rescale: when the points are at distance of order $r$ from $\g$, the centers of the auxiliary circles involved are at distance of order $r^2$. See Figure \ref{triinfty} again. 

It turn out that these orbits at infinity of the circle centers have the same shape as the orbits of points at infinity of the outer (area) billiard on the same oval $\g$. We describe these orbits and explain this somewhat unexpected phenomenon in Section \ref{sect:circinfty}. Unlike Section \ref{sect:dyninfty} where we provide rigorous proofs of our results, in Section \ref{sect:circinfty} we restrict ourself to an informal discussion: the proof would involve tedious estimates, similar to the ones in Section \ref{sect:dyninfty}, and keeping the size of the paper under control, we decided to omit them. For similar results concerning the polygonal tables, see \cite{EC}.

\bigskip

{\bf Acknowledgements}. PA acknowledges funding by the Deutsche
Forschungsgemeinschaft (DFG, German Research Foundation) through Germany's
Excellence Strategy EXC-2181/1 - 390900948 (the Heidelberg STRUCTURES Excellence
Cluster), the Transregional Colloborative Research Center CRC/TRR 191
(281071066).  LEC and ST were supported by NSF grant DMS-2404535 and Simons
Foundation grant MPS-TSM-00007747. ST was also supported by a Mercator
fellowship within the CRC/TRR 191. He and LEC thank the Heidelberg University
for its hospitality and the inspiring atmosphere.

\section{Generating function and invariant area form} \label{sect:def}

Let $\g(x)$ be arc length parameterization of the oval.
Recall from \cite{AT1} that a generating function for the map $F$ is $H(x,x')=|xA|+|Ax'|$:  if points $x$ and $x''$ are fixed and $x'\in\g$ varies, then $x'$ is a critical point of the function $H(x,x')+H(x',x'')$ if and only if $x'$ is the tangency point of a circle that is tangent to the positive tangent ray to $\g$ at $x$, the negative tangent ray at $x''$, and $\g$, as in Figure \ref{fourt}. In formulas,
$F(x,x')=(x',x'')$ if and only if
\begin{equation} \label{eq:gen}
\frac{\partial (H(x,x')+H(x',x''))}{\partial x'} =0.
\end{equation}

Since the results of the next two lemmas were stated in \cite{AT1} without proofs,  we present proofs here.

\begin{lemma} \label{lm:twist}
The partial derivatives if the generating function are as follows:
\begin{equation} \label{eq:twist}
\begin{split}
\frac{\partial H(x,x')}{\partial x}= -&k(x)|Ax| \cot \frac{\varphi}{2}-1, \ 
\frac{\partial H(x,x')}{\partial x'}= k(x')|Ax'| \cot \frac{\varphi}{2}+1, \\
&\frac{\partial^2 H(x,x')}{\partial x\partial x'}= - \frac{k(x)k(x')(|Ax|+|Ax'|)}{2\sin^2 \frac{\varphi}{2}},
\end{split}
\end{equation}
where $k$ is the curvature function of $\g$, and $\varphi$ is the angle between the tangent segments $Ax$ and $Ax'$.
\end{lemma}

\begin{figure}[ht]
\centering
\includegraphics[width=.35\textwidth]{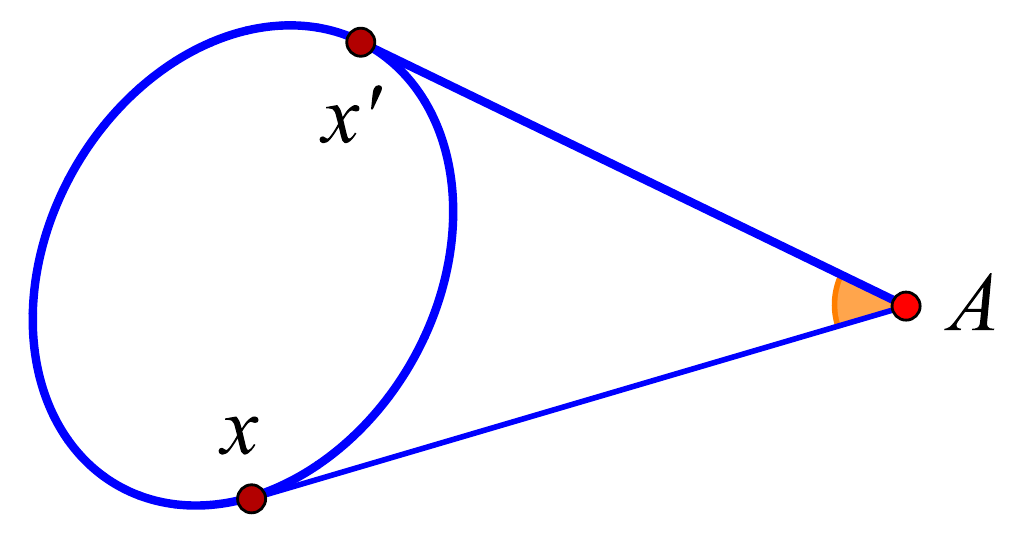}
\caption{To Lemma \ref{lm:twist}.}	
\label{der}
\end{figure}

\proof Referring to Figure \ref{der}, choose the coordinate system so that point $x'$ is the origin and $x'A$ is the positive horizontal coordinate axis. Assume that the curve near point $x$ is arc length parameterized and its components are functions $a(x)$ and $b(x)$. Then
$$
a'=\cos \phi, b'=\sin \phi, a''=-k(x) \sin \phi, b''=k(x) \cos \phi,
$$
where $k$  is the curvature and prime means $d/dx$.

Let $|xA|=\lambda$. We have
$$
a+\lambda a' = A, \ b+\lambda b' = 0.
$$
Differentiating the second equation, we get
$$
1+ \lambda' = -\frac{\lambda k \cos \phi}{\sin \phi}, 
$$
hence  
$$
\lambda' = -\frac{\lambda k \cos \phi}{\sin \phi} -1.
$$
Differentiating the first equation, yields
$$
A' =  (1+ \lambda')\cos\phi - \lambda k \sin \phi = -\frac{\lambda k}{\sin \phi}.
$$
Hence
$$
\frac{dH}{dx} = - |Ax| k(x) \cot (\phi/2) -1,\ \frac{dH}{dx'} = |Ax'| k(x') \cot (\phi/2) + 1,
$$
where the second equation is obtained the same way as the first one. 

Differentiating the second equation with respect to $x$, using the fact that $d\varphi/dx = k(x)$ and some trigonometry, yields the last formula in (\ref{eq:twist}).
\proofend

Equation (\ref{eq:twist}) implies that 
$$
\frac{\partial^2 H(x,x')}{\partial x\partial x'} < 0,
$$
the monotone twist condition for the map $F$. The map $F$ preserves the area 2-form 
$$
\omega=-\frac{\partial^2 H(x,x')}{\partial x \partial x'}\ dx\wedge dx'.
$$

\begin{remark}
{\rm Along with the arc length parameterization of the curve $\gamma$, one can parameterize it by the direction of the outward normal. Let $\alpha$ and $\alpha'$ be these directions at points $\gamma(x)$ and $\gamma(x')$. 
Then $0<\alpha' - \alpha < \pi$. 
Since $k(x) = d\alpha/dx$ and likewise for $x'$, and since $\varphi=\pi - (\alpha'-\alpha)$,
the last formula in (\ref{eq:twist}) implies:
$$
\omega=- \frac{(|Ax|+|Ax'|)}{2\cos^2 \left(\frac{\alpha'-\alpha}{2}\right)}\ d\alpha \wedge d\alpha'.
$$ 
}
\end{remark}

The form $\omega$  is a multiple of the standard area form $\omega_0$.

\begin{lemma} \label{lm:form}
One has
$$
\omega_A=\cot \left(\frac{\varphi}{2}  \right) \left( \frac{1}{|Ax|}+\frac{1}{|Ax'|} \right) \omega_0.
$$
\end{lemma}

\proof
This follows from Crofton's formula, see formula (2.8) in \cite{Sa}:
$$
\omega_0 = \frac{|Ax||Ax'| k(x) k(x')}{\sin\varphi} dx\wedge dx'
$$
and the last formula in (\ref{eq:twist}).
\proofend

\begin{remark} 
{\rm The invariant area forms of the conventional, inner, and of the outer billiard maps are independent of the shape of the billiard table: in the former case, this is the canonical symplectic form on the space of oriented lines, and in the latter it is the standard area form in the plane. Thus in these  cases the area form and the map are decoupled. 

For the conventional billiards this fact implies that the billiard maps with respect to confocal ellipses commute, see, e.g., \cite{Ta}. This is a non-trivial configuration theorem, see Figure \ref{commute}. The idea of the proof is that a confocal family of ellipses define a foliation of the space of oriented lines, inducing the Arnold-Liouville coordinates on the leaves. The billiard maps are parallel translations in these coordinates, and parallel translations commute.  

\begin{figure}[ht]
\centering
\includegraphics[width=.35\textwidth]{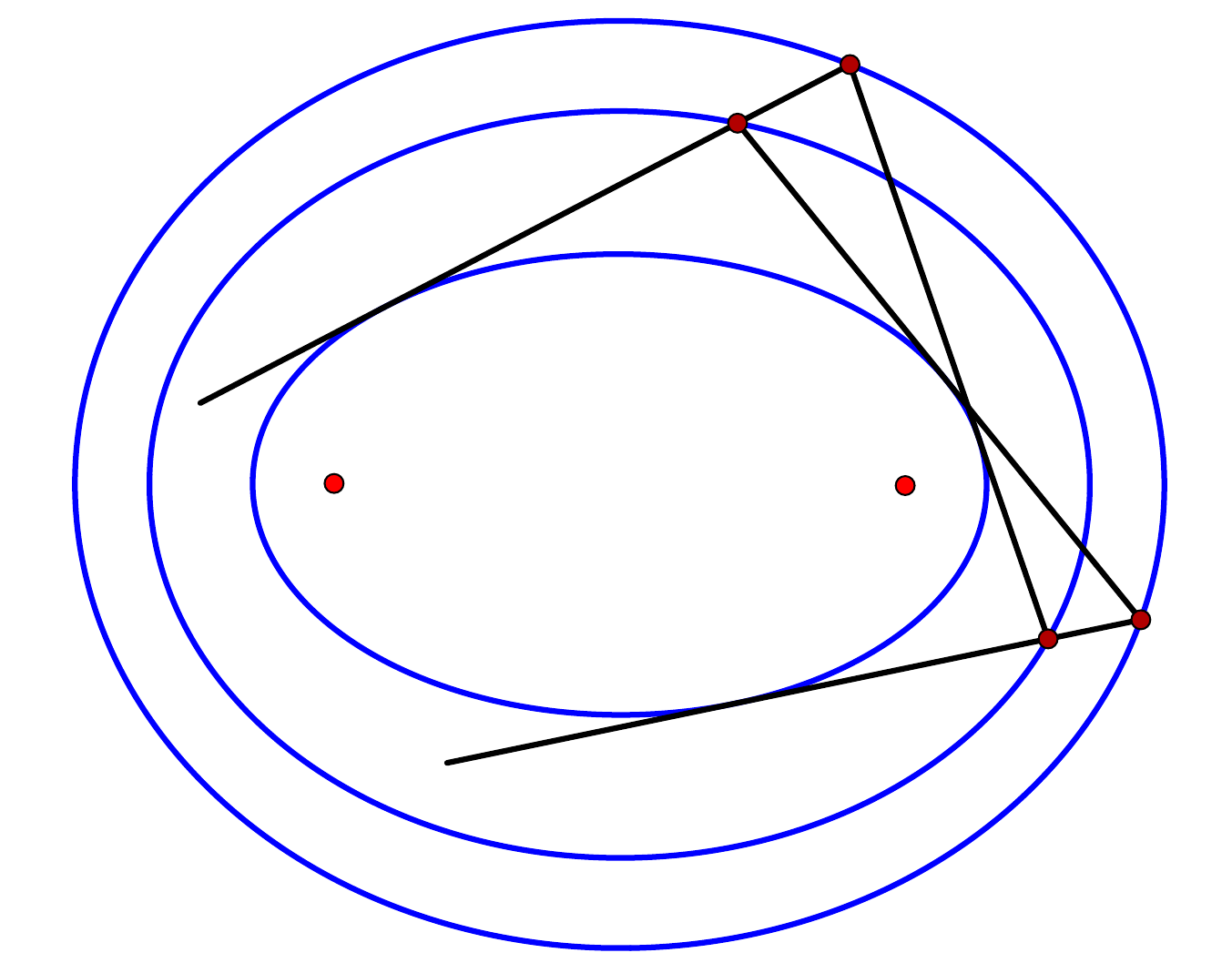}\qquad\quad
\includegraphics[width=.3\textwidth]{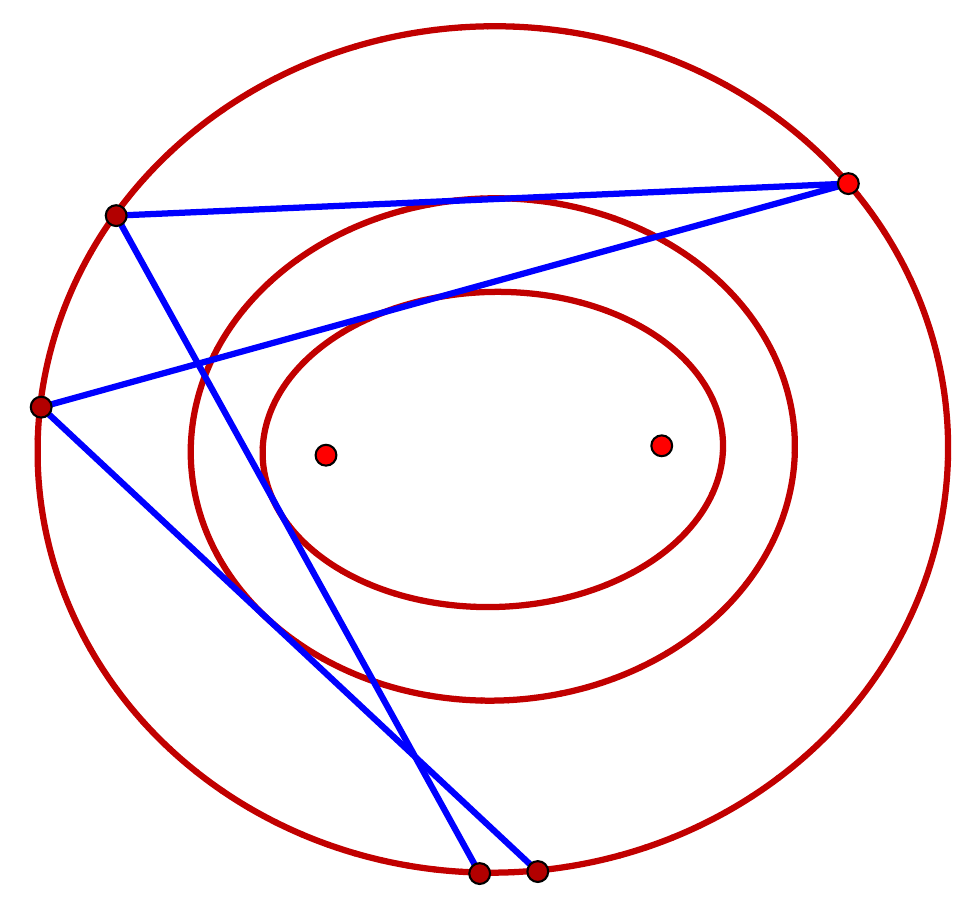}
\caption{Left: the billiard maps with respect to confocal ellipses commute. Right: the outer length billiard maps about confocal ellipses do not commute.}	
\label{commute}
\end{figure}

In contrast, for outer length billiards, as well as for symplectic billiards (see \cite{AT}), the invariant form depends on the billiard table. Although the outer length billiard map about an ellipse is integrable -- the invariant curves again are confocal ellipses (see the next section) -- the respective maps do not commute.}
\end{remark}

Next we determine the behavior of the area form $\omega$ at infinity. 
To do that we work in polar coordinates and let the points go to infinity radially. In the notations of Lemma \ref{lm:form}, we write $A=(r,\psi)$ and we consider the limit of the factor 
$$
C_A (\psi) := \cot \left(\frac{\varphi}{2}  \right) \left( \frac{1}{|Ax|}+\frac{1}{|Ax'|} \right)
$$
as $r \to \infty$. 

Let $\g$ be an oval and $\psi$ a direction. There are two  support lines of $\g$ parallel to $\psi$. The distance between them is the {\it width} $w_\gamma(\psi)$.

\begin{lemma} \label{lm:fforminfty}
One has 
 $$
 \lim_{r\to\infty} C_A(\psi) = \frac{4}{w_\gamma(\psi)}.
 $$ 
\end{lemma}

\begin{figure}[ht]
\centering
\includegraphics[width=.4\textwidth]{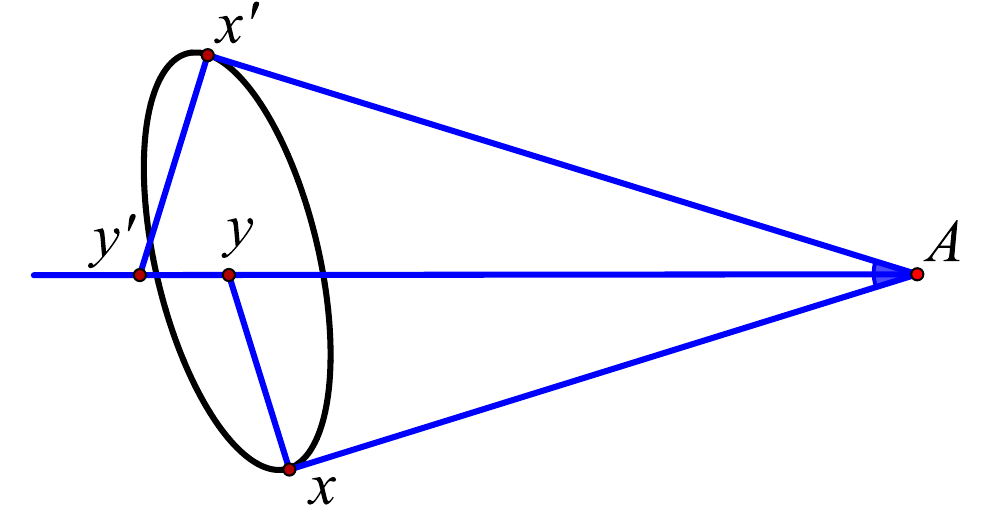} \qquad\qquad
\includegraphics[width=.4\textwidth]{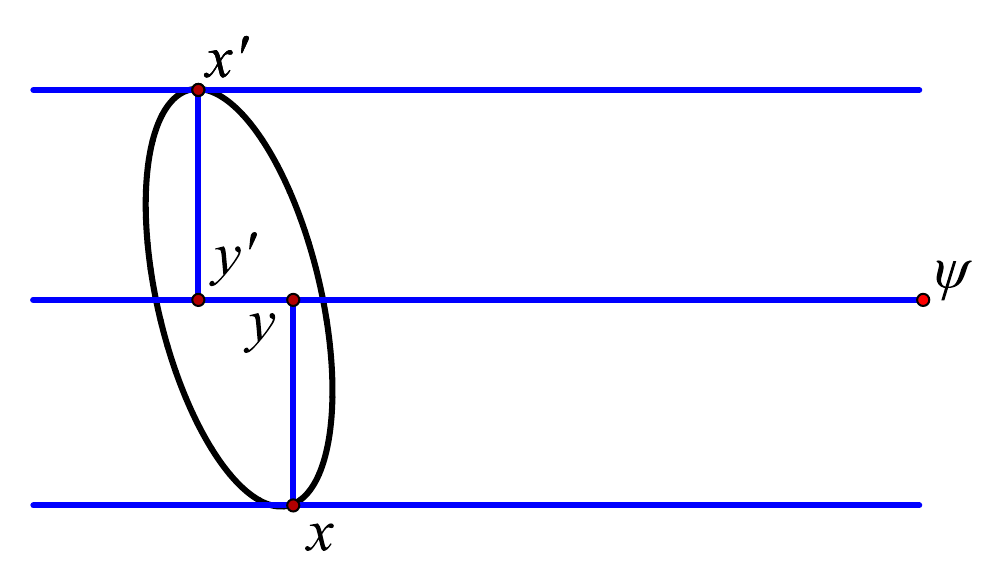}
\caption{To Lemma \ref{lm:fforminfty}.}	
\label{lim}
\end{figure}

\proof
Consider Figure \ref{lim} left. The segments $xy$ and $x'y'$ are orthogonal to the tangents $Ax$ and $Ax'$, the points $y$ and $y'$ lie on the bisector of the angle $xAx'$. Hence 
$$
\cot \left(\frac{\varphi}{2}  \right)  \frac{1}{|Ax|} = \frac{1}{|xy|},\ \ 
\cot \left(\frac{\varphi}{2}  \right) \frac{1}{|Ax'|} = \frac{1}{|x'y'|}.
$$

As $r \to \infty$, the lines $Ax$ and $Ax'$ become parallel and the bisector becomes the middle line of the strip made by these lines, see Figure \ref{lim} right. In this limit, one has $|xy|=|x'y'|=\frac{w_\gamma(\psi)}{2}$, and the result follows.
\proofend

\begin{remark}
{\rm One can analyze the limiting behavior of the form $\omega_A$ as the point tends to the curve $\g$, that is, when points $x$ and $x'$ merge together. In this limit, the factor 
$$
\cot \left(\frac{\varphi}{2}  \right) \left( \frac{1}{|Ax|}+\frac{1}{|Ax'|} \right)
$$
equals $2 k(x)$. Since we do not use this result, we omit the proof.
}
\end{remark}

\section{Ellipses} \label{sect:ell} 

Similarly to the conventional billiard, the outer length billiard map $F$ is integrable if the table $\g$ is an ellipse. The next result was proved in \cite{BBF}.

\begin{figure}[ht]
\centering
\includegraphics[width=.5\textwidth]{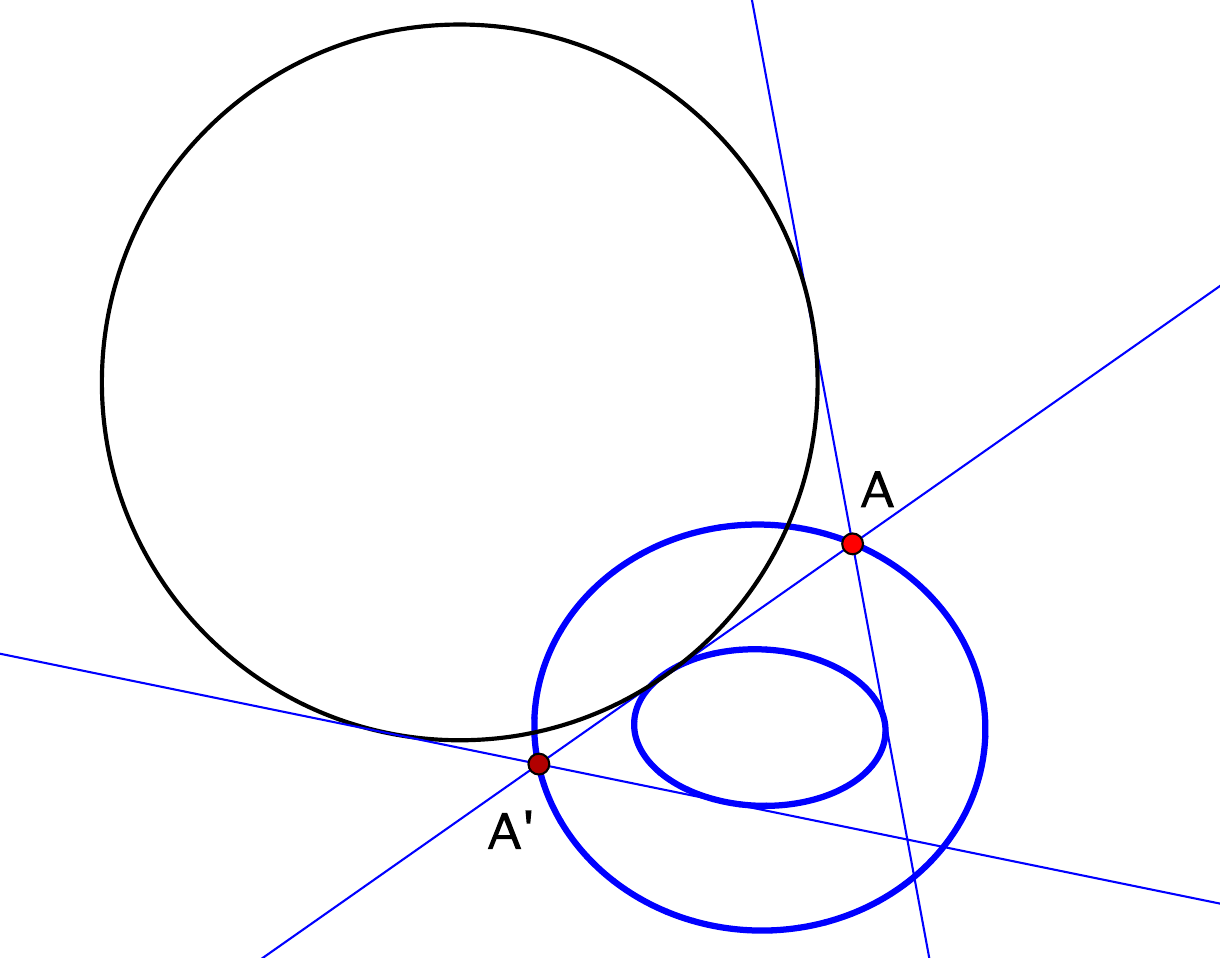}
\caption{Confocal ellipses are invariant curves of the outer length billiard about an ellipse.}	
\label{integrability}
\end{figure}

\begin{theorem} \label{thm:int}
If $\g$ is an ellipse then the confocal ellipses are invariant curves of the outer length billiard map $F$, see Figure \ref{integrability}.
\end{theorem}

In this section, we put this result into historical context and make some related comments. 
\medskip 

\begin{figure}[ht]
\centering
\includegraphics[width=.4\textwidth]{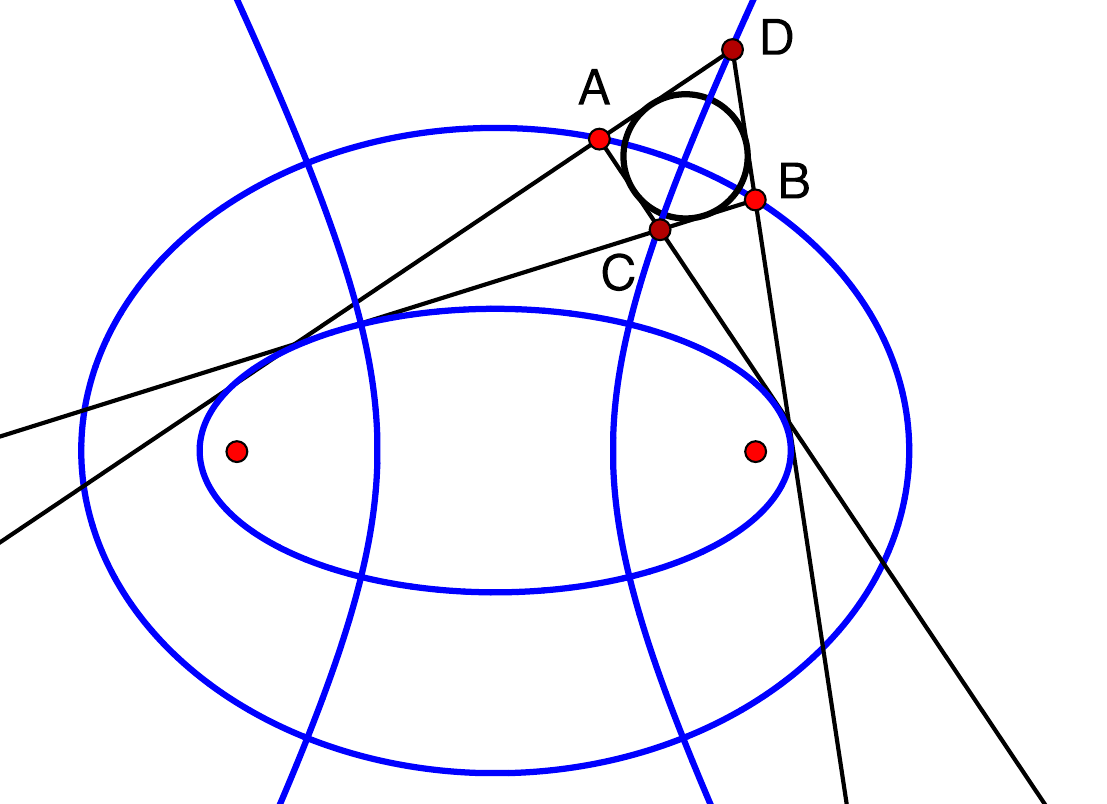}
\includegraphics[width=.59\textwidth]{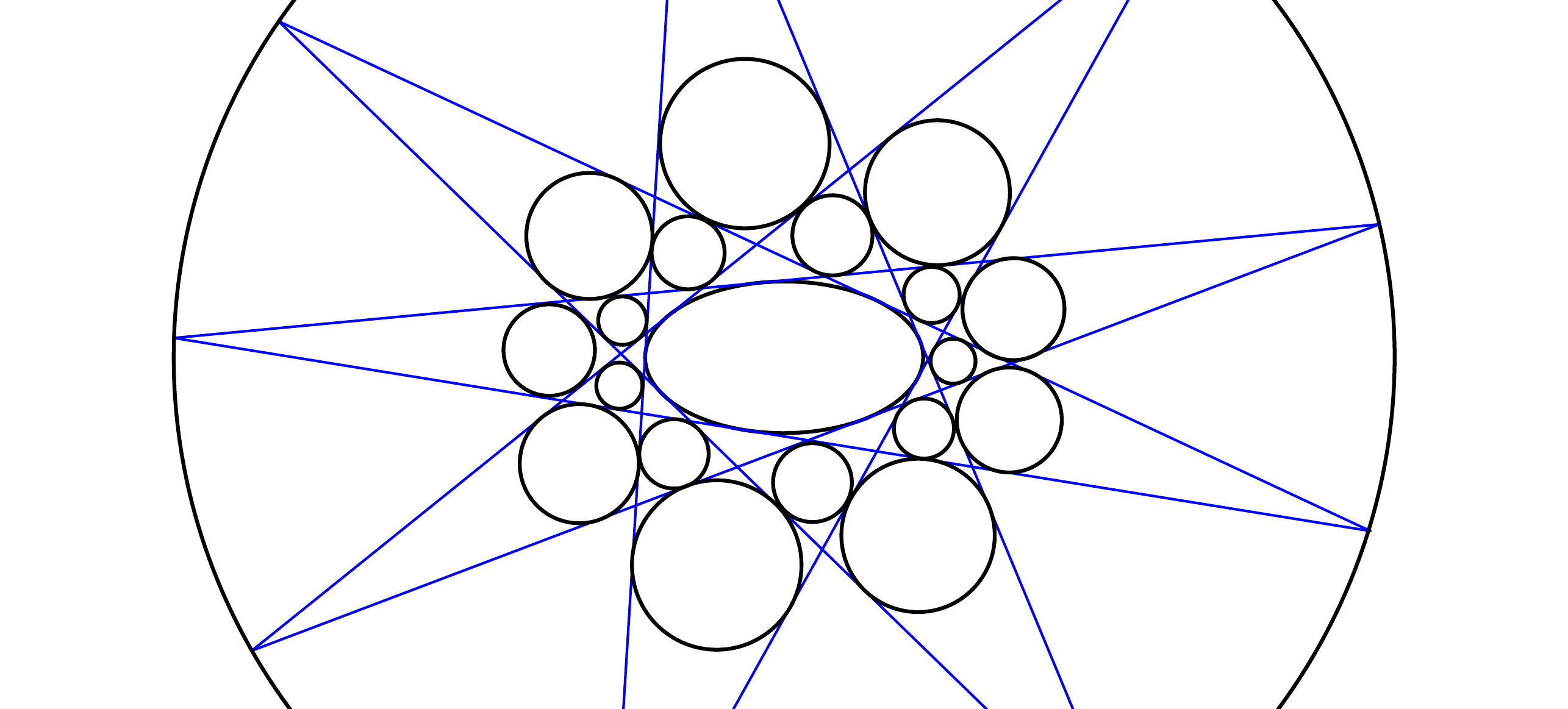}
\caption{Chasles-Reye theorem and the Poncelet grid.}	
\label{quad}
\end{figure}

(1) Figure \ref{integrability} is a limiting case of the Chasles-Reye theorem depicted in Figure \ref{quad}. Let us recall this theorem.

Let $\g$ be an ellipse and $C$ and $D$ be two points on a confocal hyperbola. Draw tangent lines to $\g$ from  points $C$ and $D$. These four lines form a quadrilateral $ACBD$. The claim of the theorem is that points $A$ and $B$ lie on a confocal ellipse, and the quadrilateral $ACBD$ is circumscribed about a circle. See \cite{AB,IT} for a modern take on this result, its generalizations, and the relation to billiards in ellipses.

If point $C$ is chosen to be on $\g$, the quadrilateral $ACBD$ degenerates to a triangle, and the circle that is tangent to the sides of this triangle is tangent to the side $AB$ at point $C$, as needed in the definition of the outer length billiard.

(2) Let $E_1 \subset E_2$ be nested confocal ellipses. Assume that there exists 
a polygon inscribed into $E_2$ and circumscribed about $E_1$, that is,  a periodic orbit of  the conventional (inner) billiard in $E_2$ and, at the same time, a periodic orbit of the outer length billiard about $E_1$. 

The extensions of the sides of such a Poncelet polygon form the Poncelet grid, see \cite{LT,Sc}. The Chasles-Reye theorem implies that the quadrilaterals of this grid are circumscribed, see Figure \ref{quad}.
\medskip

(3) The observation that the vertices of the polygons of minimal perimeter circumscribed about an ellipse lie on a confocal ellipse is due to Steiner (Werke II. pp. 618--20).\footnote{To quote M. Berry: ``Nothing is ever discovered for the first time", see \url{https://michaelberryphysics.wordpress.com/quotations/}
}
\medskip

(4) Consider again nested confocal ellipses $E_1 \subset E_2$. Then $E_2$ is an invariant curve of the outer length billiard map $F$ with respect to $E_1$. Since $F$ is completely integrable,  the Arnold-Liouville theorem implies the existence of an $F$-invariant measure on $E_2$. The description of this measure goes back to Bertrand's proof of the Poncelet porism ;  see, e.g., \cite{Av} for a recent treatment.

We present an explicit formula, a particular case of Theorem 1.5 in \cite{Av} when the ellipses are confocal. Let
$$
\frac{x^2}{a^2} + \frac{y^2}{b^2}=1\ \ \  {\rm and}\ \ \ \frac{x^2}{a^2+\lambda} + \frac{y^2}{b^2+\lambda}=1,\ \lambda >0
$$
be the equations of $E_1$ and $E_2$ respectively, and let $s$ be the arc length parameter on $E_2$. Then the $F$-invariant measure $\mu$ on $E_2$ is given by $ds/G(x,y)$ where
$$
G(x,y)={\sqrt{
\left[\frac{x^2}{(a^2+\lambda)^2}+ \frac{y^2}{(b^2+\lambda)^2}\right]
\left[\frac{x^2}{a^2} + \frac{y^2}{b^2}-1\right]}}.
$$

In particular, consider the case of  $\lambda \to \infty$. In this limit, $E_2$ is  a circle parameterized as $\sqrt{\lambda}(\cos\alpha, \sin\alpha)$. The (rescaled) invariant measure is then given by
$$
\mu=\frac{d\alpha}{\sqrt{b^2\cos^2\alpha + a^2 \sin^2\alpha}}.
$$
We note that the denominator is the support function of the ellipse $E_1$, 
see, e.g., Lemma 3.1 of \cite{BT}. This is half the width of $E_1$ in the direction $\alpha$, that is, 
$$
\mu = 2 \frac{d\alpha}{w_{E_1}(\alpha)}.
$$

Let the ellipse degenerate to a segment of length $2a$, that is, let $b=0$ (cf. Figure \ref{segment}). Then 
\begin{equation*} \label{eq:meas}
\mu=\frac{d\alpha}{a\sin\alpha},
\end{equation*}
where the denominator is half the length of the projection of the segment along the direction $\alpha$.
This measure has singularities at $\alpha=0,\pi$, corresponding to the fact that the billiard trajectory in an ellipse that passes through its foci converges to the major axis, traversed back and forth.
The vector field dual to $\mu$, namely, $a \sin\alpha\ \partial \alpha$, will play a  role in the next section.

\section{The second iteration of the outer length billiard map at infinity} \label{sect:dyninfty}

As stated in the introduction, the second iteration of the outer length billiard map at infinity  is approximated by the flow of a Hamiltonian vector field whose Hamiltonian function is the distance from the origin with respect to the area form, invariant under the map. In this section we explain and quantify this phenomenon.

Let $(\alpha,r)$ be the polar coordinates; since $\g$ is fixed, we shorthand $w_\g(\alpha)$ to $w(\alpha)$. The standard area form is $r dr \wedge d\alpha$. Consider the area form 
$$
\Omega=\frac{r dr \wedge d\alpha}{w(\alpha)}.
$$
This is, up to a factor, the invariant area form from Lemma \ref{lm:fforminfty}, extended to the whole plane. 

Consider a Hamiltonian function $H(\alpha,r)=2r$, and let 
$$
X=-\frac{2w(\alpha)}{r} \frac{\partial}{\partial \alpha}
$$
be its Hamiltonian vector field. This field is homogeneous of degree zero and its trajectories are the origin centered circles. Denote by $\Phi$ the time-$1$ flow of $X$.

As stated above, the flow of $X$ approximates the second iteration of the outer length billiard map at infinity. More precisely, one has

\begin{theorem} \label{thm:main1}
There is a constant $C$, depending only on the oval $\g$, such that 
$$
|F^2(x)-\Phi(x)| \leq \frac{C}{|x|}.
$$
\end{theorem}

The aim of this section is to prove this theorem. Before immersing in rather tedious  estimates, we describe a ``blueprint" of the proof.

\begin{figure}[ht]
\centering
\includegraphics[width=.75\textwidth]{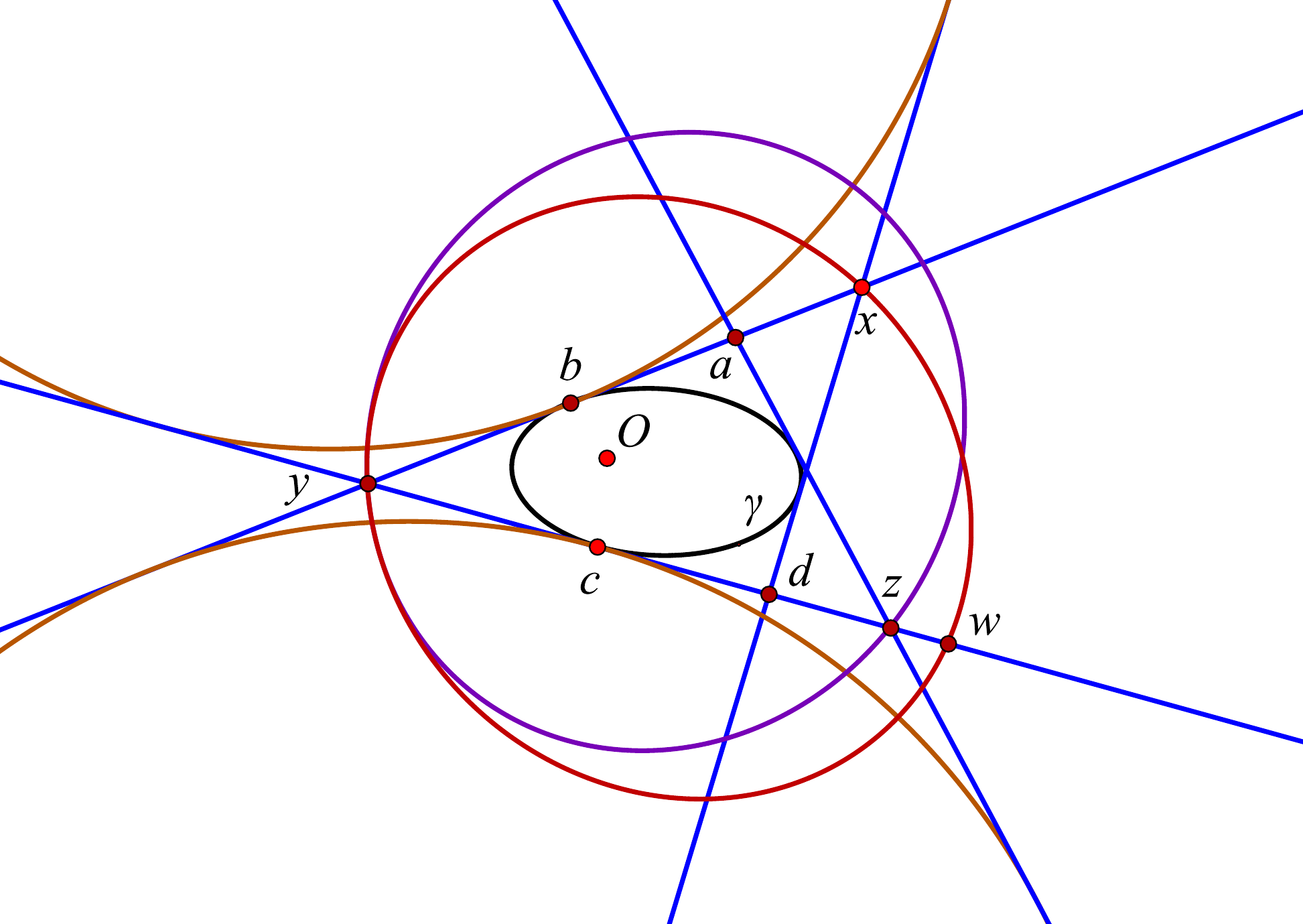}
\caption{Two ellipses.}	
\label{fig:twoell}
\end{figure}

Consider Figure \ref{fig:twoell} where $y=F(x), z=F(y)$. Since the two tangents to a (brown) circle from the same point have equal lengths, we have $|xb|+|xd|=|yb|+|yd|$, that is, points $x$ and $y$ lie on the (red) ellipse with foci $b$ and $d$. Likewise, points $y$ and $z$ lie on the (purple) ellipse with foci $a$ and $c$. 

What happens when point $x$ goes to infinity in  direction $\alpha$?
The lines $xy$ and $zy$ become nearly parallel. Points $a$ and $b$ merge together, and so do $c$ and $d$; the two resulting ``double points", say, $p$ and $q$, are the tangency points of the support lines of $\g$ in the direction $\alpha$. The two ellipses also merge together; they have  foci $p$ and $q$. Thus the outer length billiard map on the oval $\g$  is approximated by that on the segment $pq$. 

Furthermore, the large ellipses with fixed foci become nearly circular. Thus, in the limit, we obtain a circle with a segment $pq$ inside it and a Poncelet-like circle map using, alternating, the support lines of the segment at its end points. This map can be explicitly analyzed: it is approximated by the time-1 flow of the vector field that is the limit of the field $X$ as $\g$ tends to a segment. 
\medskip

Now we get to the estimates needed for the proof of Theorem \ref{thm:main1}. The constants involved  depend only on the oval $\g$ and (in some case) on the choice of the origin inside it. 

To keep the size of the paper reasonable, occasionally we use the big O notation. 
Recall the definition. Given two functions, $f(x)$ and $g(x)$, one writes
 $$
 f(x) = O(g(x))\ \ {\rm as}\ \ x\to\infty
 $$
 if there exist constants $R>0$ and $C$ such that $|f(x)| \le C|g(x)|$ for all $x$ satisfying  $|x| \ge R$. Likewise, 
 $$
 f(x) = O(g(x))\ \ {\rm as}\ \ x\to 0
 $$
  if there exist constants $\delta>0$ and $C$ such that $|f(x)| \le C|g(x)|$ for all $x$ satisfying  $|x| \le \delta$. 
\smallskip

The Hamiltonian vector field $X$ depends on the choice of the origin, and we start by describing its dependence on this choice. The next lemma makes it possible to choose the origin in a convenient way in what follows. 

Let $D$ be the diameter of the ``billiard table" $\gamma$. Let $O$ and $O'$ be two choices of the origin inside $\g$; denote the respective vector fields by $X$ and $X'$. 

\begin{lemma} \label{lm:orig}
There is a constant $C_1$ such that 
$$
|X(x)-X'(x)| \le \frac{C_1}{|x|}.
$$
\end{lemma}

\begin{figure}[ht]
\centering
\includegraphics[width=.4\textwidth]{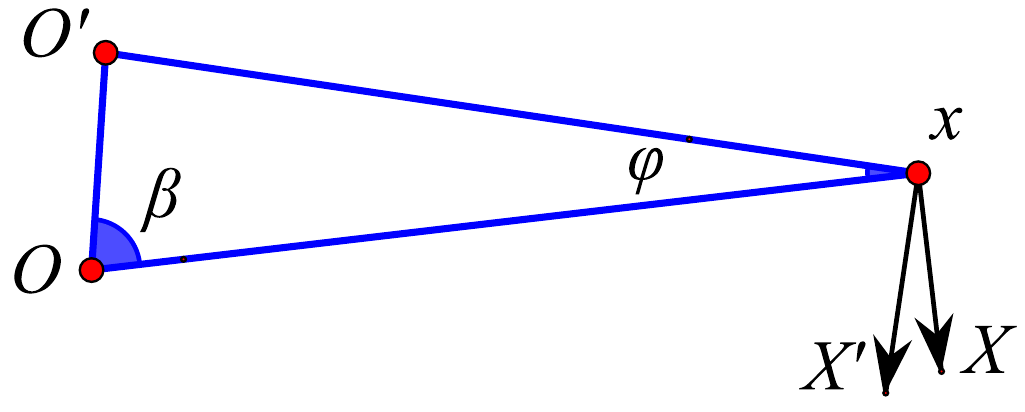}
\caption{To Lemma \ref{lm:orig}.}	
\label{fig:orig}
\end{figure}

\proof
Consider Figure \ref{fig:orig}. Let $\varphi=\alpha'-\alpha$. The Sine Rule applied to  triangle $OxO'$ implies
$$
\frac{\varphi}{2D} \le \frac{\sin\varphi}{D} \le \frac{\sin\varphi}{|OO'|} = \frac{\sin\beta}{|x|} \leq \frac{1}{|x|},
$$
where the first inequality is due to the inequality $\varphi \le 2\sin\varphi$. That is, using the big O notation, 
$$
\varphi = O\left(\frac{1}{|x|}\right).
$$

This $\varphi$ is the angle between the vectors $X(x)$ and $X'(x)$. Their magnitudes are $w(\alpha)$ and $w(\alpha')$. Let $p(\alpha)$ be the support function of the oval $\g$. Then
$$
w(\alpha)=p(\alpha)+p(\alpha+\pi),\ w(\alpha')=p(\alpha+\varphi)+p(\alpha+\pi+\varphi) = w(\alpha) + O(\varphi)
$$
(we used the fact that $\g$ is compact). 
It follows that 
$$
||X(x)|-|X'(x)|| = O\left(\frac{1}{|x|}\right)
$$
as well, and this implies the statement of the lemma. 
\proofend

 We are concerned with the orbits of the map $F$ and $F^2$ far away from $\g$, and we choose $|y|$ in Figure \ref{fig:twoell} as the measure of this large distance. Set $\eps=1/|y|$.  
 
We start with two estimates involving angles and distances related to an oval. 

\begin{figure}[ht]
\centering
\includegraphics[width=.45\textwidth]{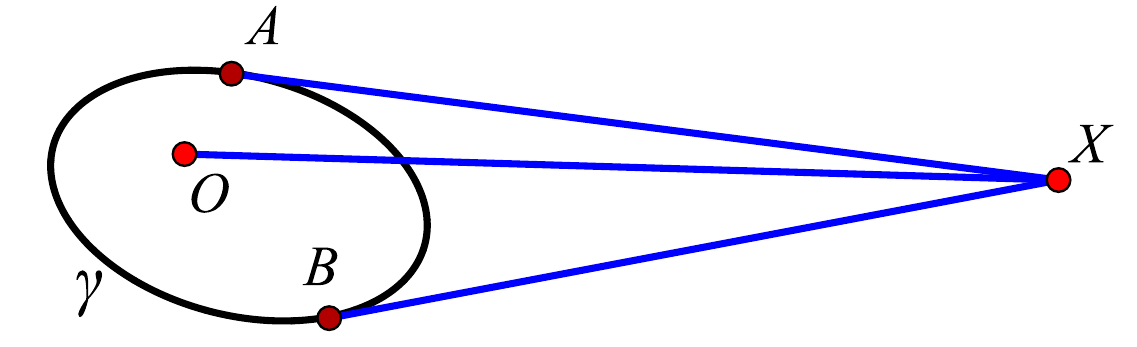}
\caption{To Lemma \ref{lm:ang1}.}	
\label{ang1}
\end{figure}

\begin{lemma} \label{lm:ang1}
Referring to Figure \ref{ang1}, one has
$$
\angle AXB \leq \frac{4D}{|X|}.
$$
\end{lemma}

\proof
Let $\alpha = \angle AXO$ and $\beta = \angle AOX$. Then
$$
\frac{\sin\alpha}{D} \le \frac{\sin\alpha}{|A|} = \frac{\sin\beta}{|X|} \le \frac{1}{|X|},
$$
where the equality is the Sine Rule for the triangle $AOX$. Thus $\sin\alpha \le \frac{D}{|X|}$. Since $\alpha \le 2\sin\alpha$, we have $\alpha \le \frac{2D}{|X|}$.
The same  applies to $\angle BXO$, implying the result.
\proofend

Referring to Figure \ref{ang1}, let $\alpha$ be the direction of the vector $OX$ and let $w(\alpha)$ be the width of $\g$ in this direction. 

\begin{lemma} \label{lm:width}
There is a constant $C_2$ such that
$$
|w(\alpha) - |AB|| \le \frac{C_2}{|x|}.
$$
\end{lemma}

\proof
Let $A_1$ and $B_1$ be the points on the curve $\g$ where the tangent lines are parallel to $OX$. The normals  to $\g$ at points $A$ and $A_1$ differ by less that $\angle AXB$. Since the curvature of $\g$ is bounded away from zero, and since $\angle AXB = O(\frac{1}{|x|})$, one has $|A_1 A| = O(\frac{1}{|x|})$. Likewise for $B$ and $B_1$,
hence 
$$
w(\alpha)=|A_1B_1|=|AB| + O\left(\frac{1}{|x|}\right),
$$
as needed.
\proofend

We need another estimate of elementary differential geometric nature.

\begin{figure}[ht]
\centering
\includegraphics[width=.23\textwidth]{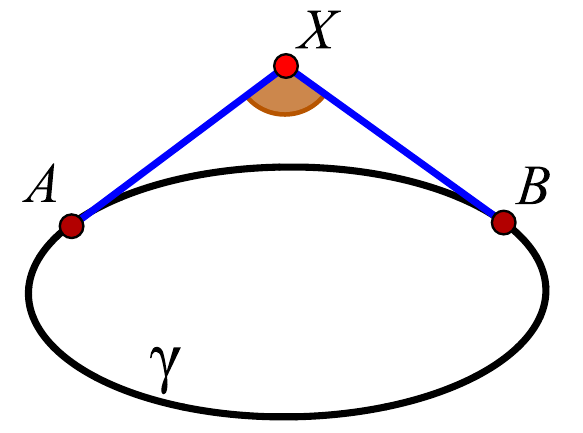}
\caption{To Lemma \ref{ang2}.}	
\label{ang2}
\end{figure}

\begin{lemma} \label{lm:ang2}
Referring to Figure \ref{ang2}, assume that $\varphi:=\angle AXB > \pi/2$. Then there exists a constant $C_3$ such that 
$$
|AX| \le C_3 (\pi-\varphi).
$$
\end{lemma}

\proof
Set $\alpha=\pi-\varphi, \beta=\angle ABX$. By the Sine Rule,
$$
\frac{\sin\beta}{|AX|} = \frac{\sin \alpha}{|AB|}.
$$
Since $\beta < \alpha < \pi/2$, one has $\sin \beta < \sin \alpha$, and hence
$$
|AX| = \frac{\sin\beta}{\sin\alpha} |AB| < |AB| < |\stackrel{\frown}{AB}|.
$$
Let $k(s)$ be the curvature of the oval $\g$ as a function of arc length, and let $m$ be the minimal value of this function. Then
$$
\alpha = \int_B^A k(s) ds \ge m |\stackrel{\frown}{AB}|.
$$
Hence 
$
|AX| \le \frac{\alpha}{m},
$
as needed.
\proofend

The next lemma shows that the tangential distances from a point and its image under the outer length  billiard map to the billiard curve differ at most by a constant.

\begin{figure}[ht]
\centering
\includegraphics[width=.75\textwidth]{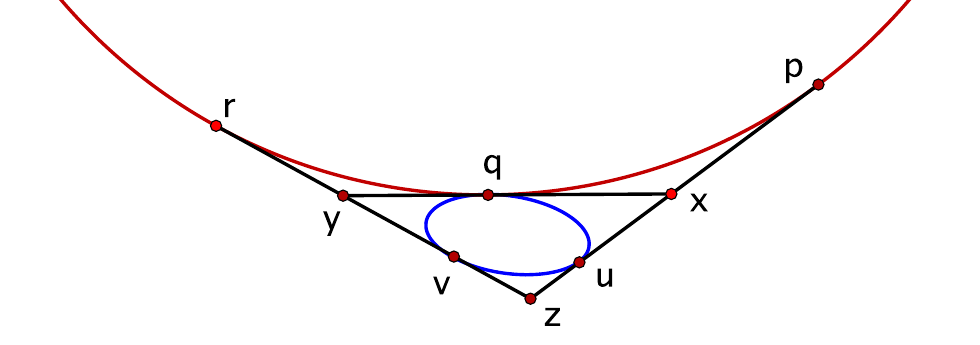}
\caption{To Lemma \ref{lm:est1}.}	
\label{fig:est1}
\end{figure}

\begin{lemma} \label{lm:est1}
Let $F(x)=y$, see Figure \ref{fig:est1}. Then
$||xq|-|yq|| \le \frac{3}{2} D$.
\end{lemma}

\proof
We use the triangle inequality and the fact that two tangents to a circle from the same point have equal lengths. 

One has
$$
|xq|-|yq|=|xp|-|yr|
$$
and
$$
|zu|+|ux|+|xp|=|zp|=|zr|=|zv|+|vy|+|yr|.
$$
Hence
$$
|xq|-|yq|=|zv|+|vy|-|zu|-|ux|,
$$
and therefore
$$
2(|xq|-|yq|)=(|xq|-|ux|)+(|zv|-|zu|)+(|vy|-|yq|).
$$

Next, 
$$
|uq|\le D,\ |vq|\le D,\ |uv|\le D,
$$
hence, by the triangle inequality, 
$$
-D \le |zv|-|zu| \le D,\ -D \le |xq|-|ux| \le D,\ -D \le |vy|-|yq| \le D.
$$
Therefore 
$$
-3D \le 2(|xq|-|yq|) \le 3D,
$$
and the result follows.
\proofend

Next we show that, for every point $x$, the distance between $x$ and $F^2(x)$ is bounded above by a constant.

\begin{figure}[ht]
\centering
\includegraphics[width=.9\textwidth]{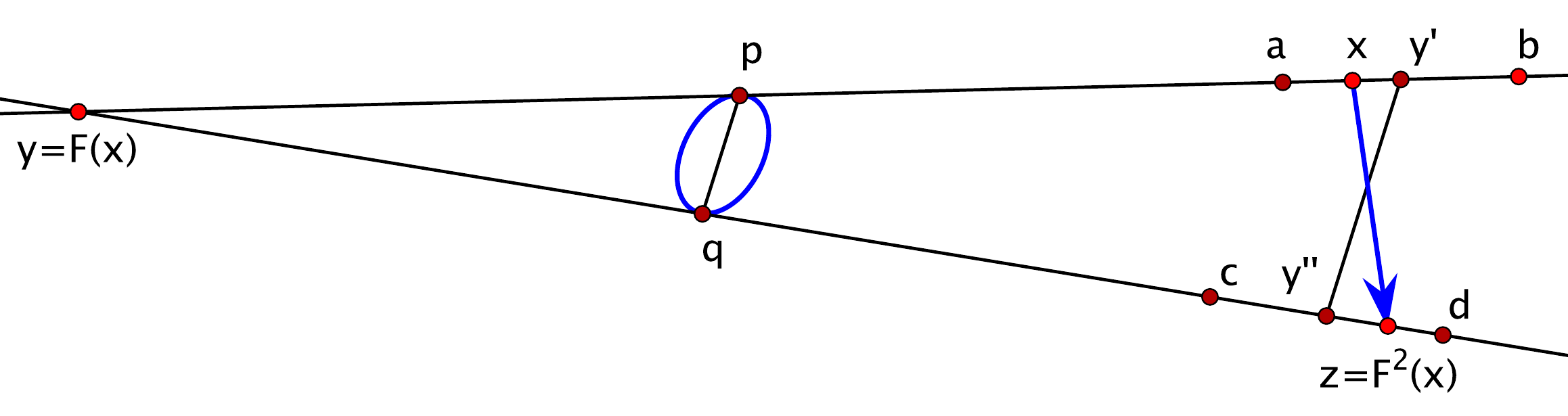}
\caption{To Lemma \ref{lm:est2}}	
\label{fig:est2}
\end{figure}

 \begin{lemma} \label{lm:est2}
Let $F^2(x)=z$, see Figure \ref{fig:est2}. Then $|xz| \le 5D$.
\end{lemma}

\proof We refer to Figure \ref{fig:est2}.
Let $y'$ and $y''$ be the reflections of point $y$ in the tangency points $p$ and $q$ respectively. Then $|y' y''|=2|pq| \le 2D$. 

Let points $a$ and $b$ be at distance $\frac{3}{2} D$ from point $y'$, and points $c$ and $d$ at distance $\frac{3}{2} D$ from point $y''$. Then, by Lemma \ref{lm:est1}, point $x$ lies on the segment $[ab]$, and point $z$ on the segment $[cd]$. By the triangle inequality,
$$
|xz| \le |xy'|+|y'y''|+|y''z|\le \frac{3}{2} D + 2D + \frac{3}{2} D = 5D,
$$
as claimed.
\proofend

The next lemma refers to Figure \ref{fig:twoell}. We use the same notations as above, and we assume that point $y$ is located so far from the origin $O$ that $|y| \ge 7D$. Recall that $\eps=1/|y|$.

\begin{lemma} \label{lm:pts}
One has: 
$$
|ab| \le 12C_3D\eps,\ |cd| \le 12C_3D\eps,
$$
where $C_3$ is the constant from Lemma \ref{lm:ang2}.
\end{lemma}

\proof
By Lemma \ref{lm:ang1}, $\angle ayd \le 4D\eps$. Using the triangle inequality, Lemma \ref{lm:est1}, and the assumption $|y| \ge 7D$, we have
\begin{equation} \label{eq:xy}
|x| \ge |xb|-|b| \ge |yb| - \frac{3}{2}D - |b| \ge |y| - \frac{3}{2}D - 2|b| \ge |y| - \frac{7}{2}D \ge \frac{1}{2}|y|,
\end{equation}
hence, by  Lemma \ref{lm:ang1}, $\angle bxd \le 8D\eps$.
Therefore 
$$
\angle ydx = \pi - \angle ayd - \angle bxd \ge \pi - 12 D \eps.
$$
By Lemma \ref{lm:ang2}, one has $|cd| \le 12C_3D\eps$. 
The same argument works for $|ab|$. 
\proofend

Next we return to Figure \ref{fig:twoell}. Let $w$ be the point of intersection of the line $yc$ with the ellipse whose foci are points $b$ and $d$ and that passes through points $x$ and $y$. 

\begin{lemma} \label{lm:zw}
One has $|zw| \le 48C_3D\eps$.
\end{lemma}

\proof
Assume that point $w$ is farther from point $y$ than point $z$ (as in the figure); the other case is treated similarly. Then since $\angle azw$ is obtuse, $|za| \le |wa|$.

For ease of notation, let $\widetilde C = 12 C_3D\eps$. Lemma \ref{lm:pts} and the triangle inequality imply
$$
|wb| \ge |wa|-\widetilde C,\ |wd| \ge |wc|-\widetilde C,\ |yd| \le |yc|+\widetilde C,\ |yb| \le |ya| +\widetilde C.
$$
Points $y$ and $z$ lie on the same (purple) ellipse with foci $a$ and $c$, and points $y$ and $w$ on the same (red) ellipse with foci $b$ and $d$, hence
$$
|yc|+|ya|=|zc|+|za|,\ |yd|+|yb|=|wd|+|wb|.
$$
Combine these equalities and inequalities:
\begin{equation*}
\begin{aligned}
&|zc|+|za| + 2\widetilde C = |yc|+|ya| + 2\widetilde C \ge |yd|+|yb|  \\
&=|wd|+|wb| \ge |wa|+|wc| - 2\widetilde C \ge |za|+|wc| -2\widetilde C,
\end{aligned}
\end{equation*}
hence 
$$
|zw|=|wc|-|zc| \le 4 \widetilde C = 48C_3D\eps,
$$
as claimed.
\proofend

Denote by $G$ the outer length billiard map on the segment $bd$ and by $H$ that on the segment $ac$. Then $$
y=F(x)=G(x),\ z=F(y)=H(y),\ w=G(y).
$$
 Therefore the above lemma implies that
$$
|F^2(x)-G^2(x)| \le 48C_3D\eps.
$$
This makes it possible to reduce the study of the second iteration of the outer length billiard map on the oval $\g$ far away from it to that on the segment $bd$: the error term is of order $\eps$.
\smallskip

Let us study the outer length billiard on a segment. 
Let $AB$ be a segment and $G$ the respective outer length billiard map. According to Section \ref{sect:ell}, the ellipses with foci $A$ and $B$ are invariant under the map $G$. Far away confocal ellipses are nearly circular, and to estimate the map $G^2$ at infinity, we approximate the invariant confocal ellipses by circles.

Normalize so that $|AB|=2$ and choose the coordinate system so that $A=(-1,0), B=(1,0)$. Then the family of confocal ellipses of given by the formula
$$
E_t = \left\{ (x,y)\big|\ \frac{x^2}{\cosh^2 t}+ \frac{y^2}{\sinh^2 t}=1 \right\}
$$
where $t>0$ is a parameter. We are interested in $t \gg 1$.

Let $C_t$ be the origin centered circle of radius $e^t/2$. 

\begin{lemma} \label{lm:Hau}
For $t\ge 2$, the Hausdorff distance between $E_t$ and $C_t$ does not exceed $e^{-t}/2$.
\end{lemma}

\proof
Parameterize $E_t$ as $(x,y)=(\cosh t \cos\alpha, \sinh t \sin \alpha)$. Then
\begin{equation*}
\begin{aligned}
\left| \sqrt{x^2+y^2}-\frac{e^t}{2}\right| = &\left| \sqrt{\sinh^2 t + \cos^2 \alpha}-\frac{e^t}{2}\right| =
\frac{\left|\sinh^2 t + \cos^2 \alpha - \frac{e^{2t}}{4}\right|}{\sqrt{\sinh^2 t + \cos^2 \alpha}+\frac{e^t}{2}}=\\
&\frac{\left| \cos^2 \alpha + \frac{e^{-2t}}{4} - \frac{1}{2} \right|}{\sqrt{\sinh^2 t + \cos^2 \alpha}+\frac{e^t}{2}}
\le \frac{  \frac{e^{-2t}}{4} + \frac{1}{2} }{\sinh t +\frac{e^t}{2}} =\frac{e^{-t}}{2},
\end{aligned}
\end{equation*}
and the result follows.
\proofend

\begin{figure}[ht]
\centering
\includegraphics[width=.55\textwidth]{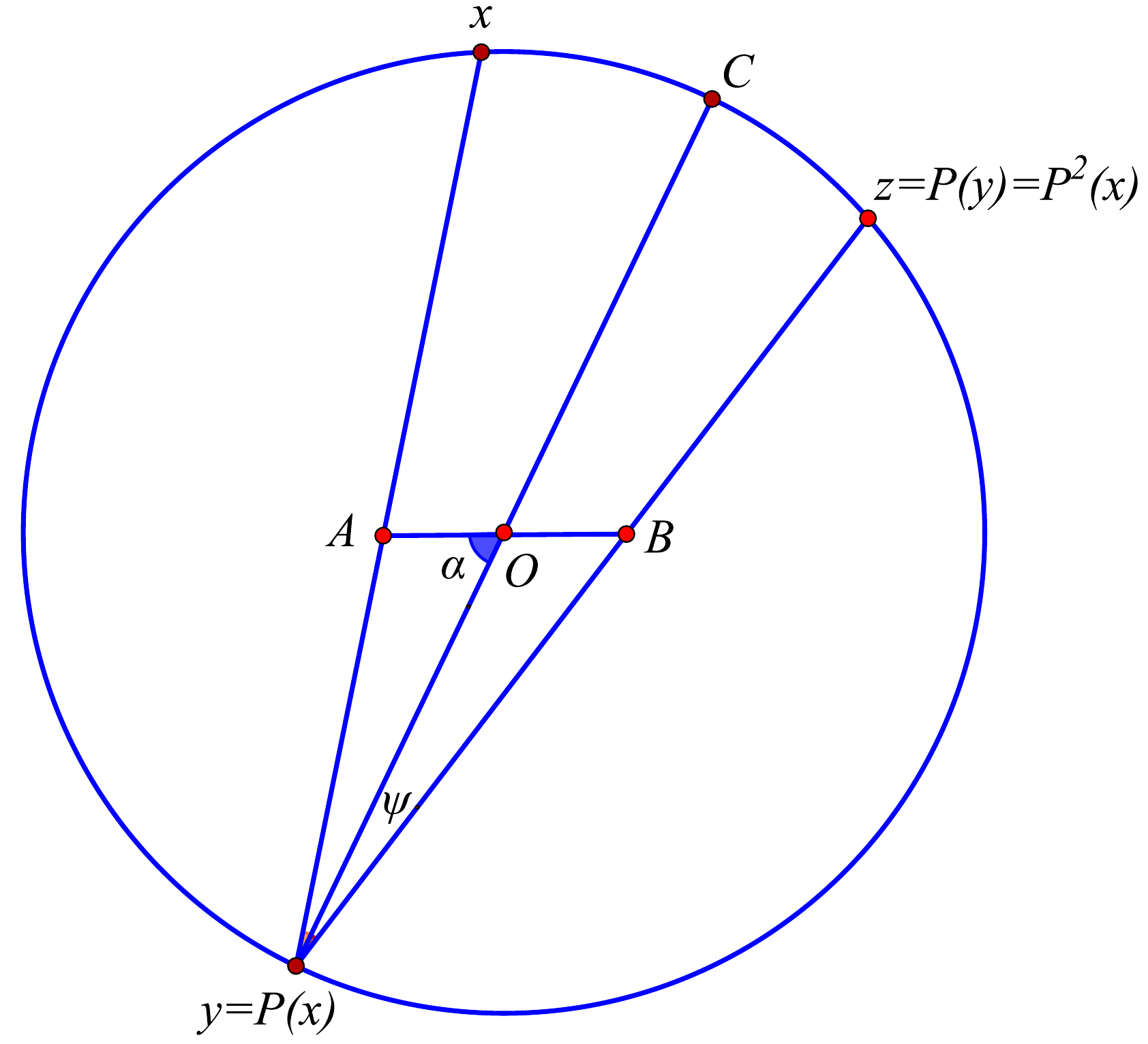}
\caption{A Poncelet map on a circle.}	
\label{Ponc}
\end{figure}

Consider  the origin centered circle of radius $r$ and denote by $P$ the ``Poncelet map"  depicted in Figure \ref{Ponc}. 

\begin{lemma} \label{lm:Ponc}
One has
$$
\angle xOz= \frac{4 \sin \alpha}{r} + O\left(\frac{1}{r^2}\right),\ \ \ {\rm and\  hence}\ \ \  |\stackrel{\frown}{xz}| = 4 \sin \alpha + O\left(\frac{1}{r}\right).
$$
\end{lemma}

\proof
The Sine Rule applied to  triangle $OyB$ yields
$$
\frac{\sin\psi}{1}= \frac{\sin(\alpha+\psi)}{r}.
$$
Hence, using some trigonometry,
$$
\sin\psi = \frac{\sin\alpha \cos\psi}{r-\cos\alpha},
$$
and therefore
$$
\psi = \frac{\sin \alpha}{r} + O\left(\frac{1}{r^2}\right).
$$
It follows that
$$
\angle COz= 2\psi = \frac{2\sin \alpha}{r} + O\left(\frac{1}{r^2}\right).
$$
Likewise for $\angle xOC$, implying the result. 
\proofend

Working in polar coordinates $(\alpha,r)$, consider the homogeneous of degree zero vector field
$$
Y=-\frac{4\sin \alpha}{r} \frac{\partial}{\partial \alpha},
$$
and let $\varphi$ be its time-$t$ flow. 

Lemma \ref{lm:Ponc} has the following corollary.

\begin{corollary} \label{cor:tan}
Referring to Figure \ref{Ponc}, one has 
$$
|P^2(x) -x - Y(x)| = O\left(\frac{1}{r}\right).
$$
\end{corollary}

Next, we estimate the displacement of a point by the map $\varphi_1$.

\begin{lemma} \label{lm:lastest}
There exist constants $C_4$ and $\delta$ such that 
$$
|\varphi_1(x)-x-Y(x)|\le \frac{C_4}{|x|}
$$
for all $x$ with $|x| \ge 1/\delta$.
\end{lemma}

\proof
The proof is similar to that of Lemma 3.3 in \cite{ACT}.

First, we claim that there exist constants $C_4\geq 0$ and $\delta>0$ such that for all $y$ with $|y|\leq 1$ and $|t|\leq \delta$ we have 
\begin{equation}\nonumber
\varphi_t(y)=y+V(y)t+R(t)\quad\text{and}\quad |R(t)|\leq C_3t^2.
\end{equation}
This is  the Taylor expansion of the map $t\mapsto\varphi_t(y)$, uniformly on the compact 1-ball. 

Next, let $S_c(x)=cx$ be the homothety with coefficient $c$. Since the vector field $Y$ is homogeneous of degree zero, one has $Y(x)=Y(cx)$. Hence its flow $\varphi_t$ satisfies
$$
\varphi_{ct} \circ S_c = S_c \circ \varphi_t \quad\text{and, in particular,}\quad \varphi_1(x)=\frac1c\varphi_{c}(cx).
$$

For $x$ with $|x|\geq 1/\delta$, set $c=\frac{1}{|x|}$ and $y=cx$. Apply the above Taylor estimate with $t=c$ to conclude that
\begin{equation}\nonumber
\begin{aligned}
&\left|\varphi_1(x)-x-Y(x)\right|=\left|\frac1c\big(\varphi_{c}(cx)-cx-cY(cx)\big)\right|
=\\
&\left|\frac1c\big(\varphi_{c}(y)-y-cY(y)\big)\right|
=\left|\frac{R(c)}{c}\right|\leq  C_4c=\frac{C_4}{|x|},
\end{aligned}
\end{equation}
as claimed.
\proofend

Combining Corollary \ref{cor:tan} and Lemma \ref{lm:lastest}, implies

\begin{corollary} \label{cor:est}
There exist constants $R$ and $C_5$ such that
$$
|P^2(x) - \varphi_1(x)| \le \frac{C_5}{|x|}.
$$
for all $x$ with $|x| \ge R$.
\end{corollary}

\paragraph{Proof of Theorem \ref{thm:main1}.} We use the notations as above.
One needs to show that 
$$
|F^2(x)-\Phi(x)| \leq \frac{C}{|x|}
$$
for a constant $C$ depending on the oval $\g$. To avoid  mentioning and combining various constants involved in the preceding lemmas, it is convenient to use the  big O notation instead. 

By Lemma \ref{lm:orig}, our estimates do not depend on the choice of the origin inside $\g$, and by the inequality (\ref{eq:xy}) from the proof of Lemma \ref{lm:pts} and a similar inequality with $x$ and $y$ swapped, the estimates involving $O(\frac{1}{|y|})$ can be replaced by ones with $O(\frac{1}{|x|})$.

We have two homogeneous of degree zero vector fields, $X$ and $Y$, and their respective time-1 flows, $\Phi$ and $\varphi_1$. By Lemma \ref{lm:width}, 
$$
X(x)=Y(x) + O\left(\frac{1}{|x|}\right).
$$
By Lemma \ref{lm:zw}, 
$$
F^2(x) = G^2(x) + O\left(\frac{1}{|x|}\right),
$$ and by Lemma \ref{lm:Hau}, 
$$
G^2(x)=P^2(x)+O\left(\frac{1}{|x|}\right).
$$
 By Corollary \ref{cor:est}, 
 $$
 P^2(x)=\varphi_1(x)+O\left(\frac{1}{|x|}\right).
 $$ 
 Lemma \ref{lm:lastest} and its analog for the vector field $X$ implies that 
$$
\varphi_1(x)-x=Y(x)+O\left(\frac{1}{|x|}\right),\  \Phi(x)-x=X(x)+O\left(\frac{1}{|x|}\right),
$$
therefore $\varphi_1(x)=\Phi(x)+O\left(\frac{1}{|x|}\right)$, concluding a long chain of equivalences mod $\frac{1}{|x|}$ and establishing the result. \proofend

\subsection{Distant periodic points have large periods} \label{subsect:per}

Next, we consider periodic orbits of the outer length billiard map $F$. What follows is similar to material related to Theorem 3 of \cite{ACT}.

Let $z$ be a $k$-periodic point of $F$. Then $z$ is a $k$-periodic point of $F^2$ as well. Assume that $|z| \gg 1$. 
Then the trajectory of the vector field $X$ through point $z$ is a circle of a large radius. Since the map $F^2$ is approximated by the time-1 flow of the vector field $X$, and since this field is homogeneous of degree 0, one expects that a large number of iterations of the map $F^2$ is needed for point $z$ to return back to the initial position. That is, the period $k$ must be large together with $|z|$.

Now we turn this informal argument into a rigorous one.

\begin{theorem} \label{thm:main2}
Given a positive integer $k$, there is a constant $\rho$ (depending on the oval $\g$) such that there are no $k$-periodic orbits of the outer length billiard map outside of the disc of radius $\rho$.
\end{theorem} 

Denote by $S$ the origin centered annulus bounded by the circles of radii $1/2$ and $3/2$. Let $m>0$ be the minimum of $|X(x)|$ in $S$ (due to the homogeneity of degree 0, $m$ is the minimum of  $|X(x)|$ in the punctured plane as well). The field $X$ is uniformly continuous in $S$ hence
$$
\exists \eta > 0\ {\rm such\  that\ if}\ x,y\in S\ {\rm satisfy}\ |x-y|<\eta,\ {\rm then}\ |X(x)-X(y)| < \frac{m}{2}.
$$
The next lemma is borrowed from \cite{ACT}, Lemma 3.5.

\begin{lemma} \label{lm:2pts}
Let $z$ and $p$ be two points in $\R^{2}\setminus\{0\}$ with $|z-p| < \eta |z|$.
Then one has $\frac{p}{|z|} \in S$ and $|X(z)-X(p)| < \frac{m}{2}$.
\end{lemma} 

\proof
We claim that
$$
\frac{1}{2} < \frac{|p|}{|z|} <\frac{3}{2}.
$$
Indeed, by the triangle inequality,
$$
|p| \le |z|+|p-z| < (1+\eta) |z| < \frac{3}{2} |z|
$$
and
$$
|p| \ge |z| - |p-z| > (1-\eta) |z| > \frac{1}{2} |z|.
$$
Therefore 
$$
\frac{z}{|z|} \in S,\  \frac{p}{|z|} \in S,\ \left|\frac{z}{|z|} - \frac{p}{|z|}\right| < \eta.
$$
Since the field $X$ is homogeneous of degree zero, 
$$
|X(z)-X(p)| = \left| X\left(\frac{z}{|z|}\right) - X\left(\frac{p}{|z|}\right)\right| < \frac{m}{2},
$$
as needed.
\proofend

\paragraph{Proof of Theorem \ref{thm:main2}.} The chain of big O equivalences in the proof of Theorem \ref{thm:main1} above imply the existence of constants $C_6$ and $R$ such that
$$
||F^2(x)|-|x|| \le \frac{C_6}{|x|}
$$
for al $x$ with $|x| \ge R$.
 Let 
$$
z:=p_1,p_2,\ldots,p_k:\ p_{i+1}=F^2(p_i),\ i=1,\ldots,k-1,\ F^2(p_k)=p_1
$$
be a $k$-periodic orbit of $F^2$. Then
\begin{equation} \label{eq:ppX}
\left| \sum_{i=1}^k X(p_i)\right| = \left| \sum_{i=1}^k p_{i+1}-p_i -X(p_i)\right|  \le
\sum_{i=1}^k  |p_{i+1}-p_i -X(p_i)| \le C_6 \sum_{i=1}^k \frac{1}{|p_i|}
\end{equation}
for $|p_i| \ge R,\ i=1,\ldots,k$. 

By Lemma \ref{lm:est2}, $|p_{i+1}-p_i| \le 5D$, and hence, by the triangle inequality, 
\begin{equation} \label{eq:zp}
|z-p_i|\le 5(k-1)D,\ i=1,\ldots,k-1. 
\end{equation}
By the triangle inequality again, this implies that if $|z| \ge R + 5(k-1)D$ then $|p_i| \ge R$ for all $i$. Thus, we 
assume that $|z| \ge R + 5(k-1)D$ to be able to  use inequality (\ref{eq:ppX}).

Recall Lemma \ref{lm:2pts} and assume that
$$
|z| > \frac{5(k-1)D}{\eta}.
$$
Combining this with \eqref{eq:zp}, we  obtain $|z-p_i| \le \eta |z|$ for all $i$. Apply Lemma \ref{lm:2pts} and conclude that $\frac{|p_i|}{|z|}\geq\frac12$, and therefore 
$$
\sum_{i=1}^k \frac{1}{|p_i|} \leq \frac{2k}{|z|}.
$$
On the other hand, Lemma \ref{lm:2pts} asserts that $|X(z)-X(p_i)| < \frac{m}{2}$ for all $i$ where $\displaystyle m=\min_{x\in\R^{2}\setminus\{0\}}|X(x)|$. In particular, $|X(z)| \ge m$, and it follows that
\begin{equation}\nonumber
\begin{aligned}
\left|\sum_{i=1}^k X(p_i)\right| = \left|\sum_{i=1}^k [X(p_i)-X(z)] +kX(z)\right|
\ge k|X(z)| - \sum_{i=1}^k |X(z)-X(p_i)|
\ge \frac{km}{2}.
\end{aligned}
\end{equation}
Thus (\ref{eq:ppX}) implies that 
$$
\frac{km}{2} \le \frac{2kC_6}{|z|}\ \ {\rm or}\ \ |z| \le \frac{4 C_6}{m}.
$$

To summarize, we showed that 
$$
{\rm if}\ \ |z| >\max\left\{R + 5(k-1)C_6,\frac{5(k-1)C_6}{\eta}\right\}\quad {\rm then} \quad |z| \leq \frac{4C_6}{m}
$$
for a $k$-periodic point $z$. Equivalently, if $z$ is a $k$-periodic point then
\begin{equation}\nonumber
|z|\leq\max \left\{R + 5(k-1)C_6,\ \frac{5(k-1)C_6}{\eta},\ \frac{4C_6}{m}\right\}=:\rho(k,\g),
\end{equation}
as needed. \proofend

\subsection{Invariant curves at infinity} 
Since the second iteration of the outer length billiard map  is approximated by the flow of a Hamiltonian vector field, one expect a KAM-type result on the existence of invariant curves ``at infinity". Indeed, assuming that $\gamma$ is sufficiently smooth and has positive curvature, we apply Moser's small twist theorem \cite{Mo} (see also \cite{La})  to establish the existence of invariant curves sufficiently far away from the billiard table. These curves prevent the orbits from escaping to infinity.

The proof is similar to the proof of the same fact for the (area) outer billiard, see, e.g., \cite{Ta95}. 
Let $(\alpha,r)$ be the polar coordinates.  Set $\rho=1/r$. Then $(\alpha,\rho)$ are the coordinates ``at infinity".

\begin{lemma} \label{lm:app1}
Let $F(\alpha,\rho)=(\alpha',\rho')$. Then
$$
\alpha'=\alpha+\pi+f(\alpha) \rho + O(\rho^2),\ \ \rho'=\rho+O(\rho^2),
$$
where $f(\alpha)<0$.
\end{lemma} 

\proof
We refer to Figure \ref{fig:twoell}, so $|x|=r, |y|=r'$. 
 Inequality (\ref{eq:xy}) from the proof of Lemma \ref{lm:pts} implies that 
\begin{equation} \label{eq:rr}
|r-r'|\le \frac{7}{2}D.
\end{equation}
Lemma \ref{lm:ang1} implies 
$$
\angle xOy = O\left(\frac{1}{r} \right) 
$$
therefore $\alpha'=\alpha+\pi+f(\alpha) \rho + O(\rho^2)$ with a negative $f(\alpha)$. 

Next, assume that $r \ge 7 D$. Then
$$
\frac{1}{r-\frac{7}{2}D} \le \frac{2}{r}.
$$
Using (\ref{eq:rr}), we obtain 
$$
\rho'-\rho=\frac{r-r'}{rr'}\le  \frac{\frac{7}{2}D}{r (r-\frac{7}{2}d)} \le 7D \rho^2,
$$
implying $\rho'=\rho+O(\rho^2)$.
\proofend

Next, we change coordinates near $\rho=0$ to be able to apply Lazutkin's version of the invariant curve theorem, see \cite{La}.

\begin{lemma} \label{lm:change}
Let $(u,v)$ be coordinates of the cylinder $S^1\times \R_+$ near the boundary $v=0$, and let a smooth map of the cylinder be given by
$$
u'=u+f(u)v+O(v^2),\ v'=v+g(u)v^2+O(v^3),
$$
where $f(u)\neq 0$. Then there exists a smooth change of coordinates
$$
x=a(u) +O(v),\ y=b(u)v+O(v^2)
$$
such that the map is given by
\begin{equation} \label{eq:Laz}
x'=x+y+O(y^2),\ y'=y+O(y^3).
\end{equation}
\end{lemma}

\proof
To satisfy $x'=x+y+O(y^2)=x+y+O(v^2)$, calculate modulo $v^2$:
$$
a(u)+cv+b(u)v=a(u')+cv'=a(u)+a'(u)f(u)v+cv,
$$
hence 
$$
a'=\frac{b}{f}\ \ {\rm and}\ \ a(u) = \int \frac{b(u)}{f(u)}\ du.
$$ 

To satisfy  $y'=y+O(y^3)=y+O(v^3)$, calculate modulo $v^3$:
$$
b(u)v + cv^2 = b(u')v'+c(v')^2=[b(u) + b'(u)f(u)v][v+g(u)v^2]+cv^2,
$$
hence 
$$
bg+b'f=0\ \ {\rm and}\ \ b(u) = exp \left(\int -\frac{g(u)}{f(u)}\ du  \right).
$$
This implies (\ref{eq:Laz}).
\proofend

Our billiard map is area preserving hence  every simple closed curve that goes around $\g$ 
intersects its image under the outer length billiard map. Therefore 
Moser's small twist theorem, applied to the map given by (\ref{eq:Laz}), implies the existence of invariant curves arbitrarily near the boundary $y=0$, that is, near $v=0$ in the original $(u,v)$ coordinates. Lemma \ref{lm:app1} makes it possible to apply this result to the  outer length billiard map (with $u=\alpha, v=\rho$), and we conclude that the map has invariant curves arbitrarily far away from the billiard table.

\section{Dynamics of centers of the auxiliary circles at infinity} \label{sect:circinfty}

The definition of the outer length billiard map $F$ involves auxiliary circles. In the preceding section we studied the dynamics of points at infinity under the map $F^2$. In this last section we describe how the respective centers of the auxiliary circles evolve under the second iteration of the map. We will restrict ourselves to an informal explanation; rigorous proofs are similar to the arguments  presented in Section \ref{sect:dyninfty}. In particular, we will freely use the big O notation.

Consider every other center of the auxiliary circle. After rescaling, their evolution appears as a continuous motion along a closed curve. Somewhat unexpectedly, these curves are homothetic to the trajectories of points under the second iteration of the (area) outer billiard map\footnote{That is, the map depicted in the middle of Figure \ref{three}.} at infinity, but the speed of motion is different from that in the case of the outer billiard map.

We start with a description of the dynamics of the outer billiard map at infinity. We follow the recent paper \cite{ACT} where this problem is investigated in the  multidimensional symplectic setting; see also \cite{Ta96,DT} for the planar case.

Recall two notions from convex geometry. 

Let $\g$ be an oval. Its {\it central symmetrization}, $\g \ominus \g$, is the Minkowski sum of $\g$ with its centrally symmetric oval. If $\g$ is given by its support function $p: S^1\to \R$, then the support function of $\g \ominus \g$ is $p(\alpha)+p(\alpha+\pi)$. For example, the central symmetrization of a curve of constant width is a circle. 

Let $\g$ be an oval containing the origin inside. Assign to  point $x\in \g$ the covector $y$ such that
${\rm Ker}\ y = T_x\g$ and normalized so that $y(x)=1$. The collection of these covectors, as $x$ traverses $\g$, is the {\it polar dual} curve $\g^*$,  an oval in the dual plane. Using Euclidean structure to identify vectors and covectors, $\g^*$ becomes a curve in the original plane.

\begin{figure}[ht]
\centering
\includegraphics[width=.21\textwidth]{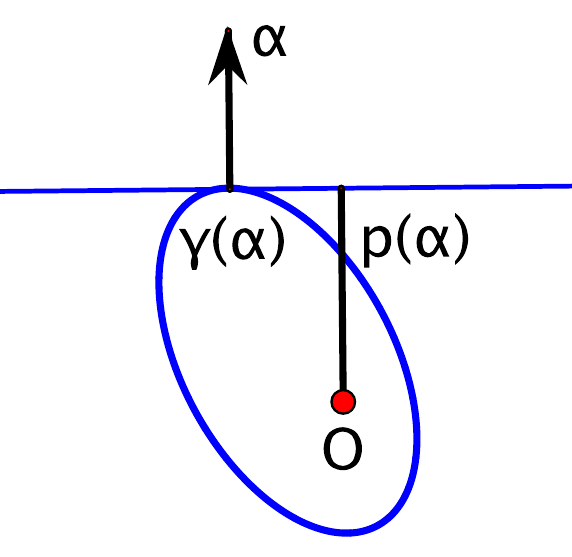}
\caption{Parameterization via the support function.}	
\label{supp}
\end{figure}

An oval can be characterized by its support function $p(\alpha)$ or by the radial function $r(\alpha)$. In the former case, one has
\begin{equation} \label{eq:supp}
\g(\alpha) = (p(\alpha)\cos\alpha - p'(\alpha \sin\alpha, p(\alpha)\sin\alpha + p'(\alpha \cos\alpha),
\end{equation}
see Figure \ref{supp}. 

In the latter case, one has 
$$
\g(\alpha) = (r(\alpha) \cos\alpha, r(\alpha) \sin\alpha),
$$
in polar coordinates.

\begin{lemma} \label{lm:par}
Let $\g$ be parameterized by its support function $p(\alpha)$ and $\g^*$ by its radial function $r(\alpha)$. Then $r(\alpha)=\frac{1}{p(\alpha)}$.
\end{lemma} 

\proof
Let $x=\g(\alpha)$, as in (\ref{eq:supp}). Then $T_x\g$ is spanned by $(-\sin\alpha,\cos\alpha)$. 
 The (co)vector $y=\frac{1}{p(\alpha)} (\cos\alpha,\sin\alpha)$ satisfies both $y^\perp = T_x\g$ and $y\cdot x=1$, as needed.
\proofend

Let $\g$ be an oval. Its width in  direction  $\alpha$ is given by $w(\alpha)=p(\alpha)+p(\alpha+\pi)$.
Applying the above lemma to the symmetrization of $\g$, we obtain

 \begin{corollary} \label{cor:dsymm}
 Let $w(\alpha)$ be the width of an oval $\g$.
 The dual to the symmetrized curve, $(\g\ominus \g)^*$, has the radial function $r(\alpha)=\frac{1}{w(\alpha)}$.
 \end{corollary}

One can also use the standard area form $\omega_0$ to identify the plane with its dual. Then the dual curve $\g^*$ again becomes a curve in the original plane, termed {\it symplectic polar dual\footnote{This construction makes sense in any symplectic space.}}: it is rotated $90^{\circ}$ compared to the polar dual obtained via the Euclidean structure.

Let $\g$ be the outer billiard oval. Consider the symplectic polar of its symmetrization, $(\g \ominus \g)^*$. Let $H$ be the homogeneous of degree 1 function whose unit level curve is $(\g \ominus \g)^*$. In view of Corollary \ref{cor:dsymm},  in polar coordinates one has
\begin{equation} \label{eq:Ham}
H(r,\alpha) = r w\left( \alpha + \frac{\pi}{2}\right).
\end{equation}
Let $X_H$ be the  Hamiltonian vector field of $H$ (with respect to $\omega_0$). 

The planar case of Theorem 1 from \cite{ACT} asserts that {\it the second iteration of the outer billiard map ``at infinity" is approximated by the flow of $X_H$} (the actual statement is similar to our Theorem \ref{thm:main1}). In particular, the trajectories of the points appear as closed curves homothetic to $(\g \ominus \g)^*$, and the motion along these curves satisfies Kepler's Second Law: the rate of change of the sectorial area is constant.

 \begin{figure}[ht]
\centering
\includegraphics[width=.7\textwidth]{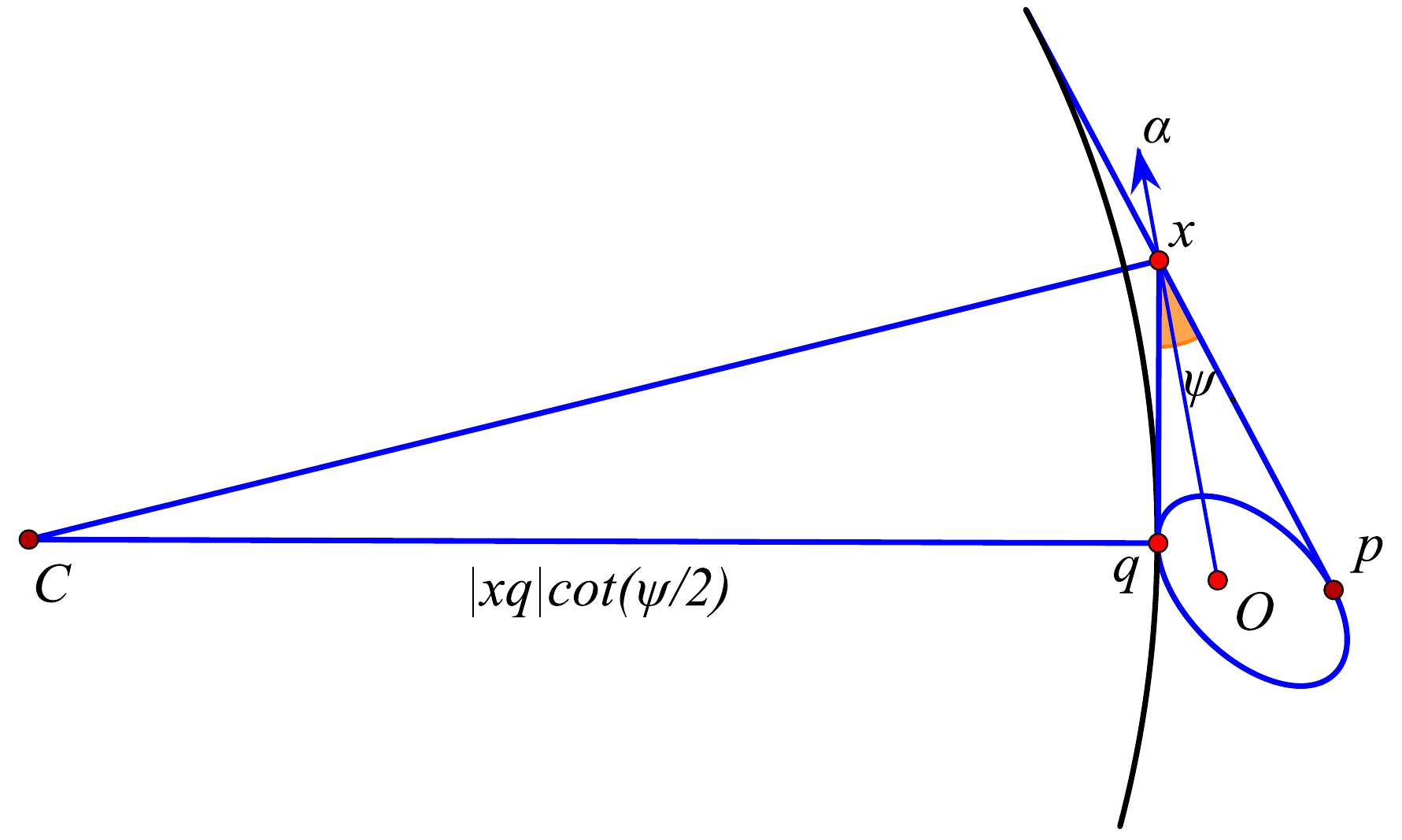}
\caption{The center $C$ of the circle involved in the definition of $F(x)$.}	
\label{center}
\end{figure}

Let us return to the outer length billiard map $F$. Consider Figure \ref{center}. Let $C(x)$ be the center of the circle involved in constructing point $F(x)$. Let $x$ and $C$ have polar coordinates $(\alpha,r)$ and $(\beta,R)$ respectively. For notational convenience, let $\eps=1/r$. 

By elementary geometry, $|qC|= |qx| \cot(\psi/2)$. 
One has
$$
|qx|=r + O(1),\ |px|=r + O(1),\ \ |qC|=R+O(1).
$$
Lemma \ref{lm:fforminfty} implies that 
$$
\cot\left(\frac{\psi}{2}\right)=\frac{2r}{w(\alpha)} +O(1),\ \ {\rm hence}\ \ R=\frac{2r^2}{w(\alpha)} +O(r).
$$
Since
$
\beta = \alpha + \frac{\pi}{2} + O(\eps),
$
one has $w(\alpha)=w(\beta+\frac{\pi}{2})+O(\eps)$. Hence 
$$
R w\left(\beta + \frac{\pi}{2}\right) = 2r^2 + O(r).
$$

Let $x'=F^2(x)$, and let $C'$ be the respective center of the auxiliary circle. We use prime to indicate that the respective quantities are related to $x'$ and $C'$.

According to Theorem \ref{thm:main1}, $r'=r+O(\eps)$, hence $(r')^2=r^2 +O(1)$. Since $\alpha'=\alpha+O(\eps)$, it follows that
$
R' w(\alpha')=R w(\alpha) + O(1).
$
Rescaling by dividing by $r^2$ (to fit into a computer screen), we conclude that, up to an error of order $\eps^2$, points $C$ and $C'$ lie on a level curve of the function (\ref{eq:Ham}), the curve traced by the orbit of a point at infinity under the second iteration of the outer billiard map. 

To find the  speed of the limiting continuous motion of the point $C$, let $t$ be the time parameter of the Hamiltonian flow that approximates the motion of point $x$ under the map $F^2$. Denote the time derivative by dot. The sectorial area rate of change  is $\frac{1}{2} R^2 \dot \beta$. In the limit $r\to\infty$,
$$
\beta = \alpha + \frac{\pi}{2},\ R = \frac{2r^2}{w(\alpha)},
$$
and since 
$$
\dot \alpha = \frac{w(\alpha)}{r},\  \dot r =0,
$$
one has 
$$
\frac{1}{2} R^2 \dot \beta = \frac{2 r^3}{w(\alpha)}.
$$
That is, the rate of change of the sectorial area is inverse proportional to the widths of the oval $\g$.

\begin{figure}[ht]
\centering
\includegraphics[width=.4\textwidth]{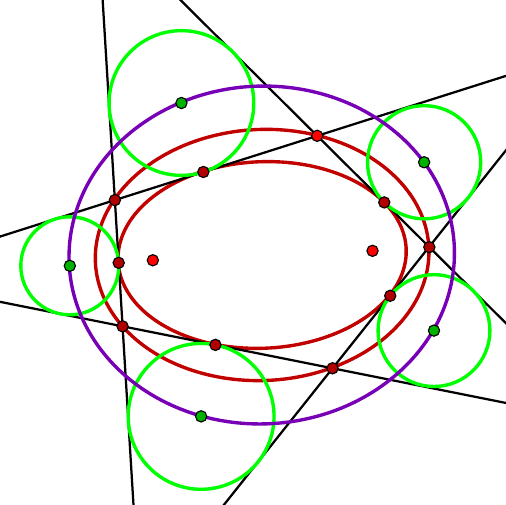}\qquad
\includegraphics[width=.35\textwidth]{Centers1}
\caption{$5$-periodic orbits of the outer length billiard map with the rotation numbers 1 and 2 and the centers of the corresponding circles. The inner ellipse is $\g$, the confocal ellipse contains the 5-periodic orbit.}	
\label{centers}
\end{figure}

\begin{example} \label{ex:ellipses}
{\rm If $\g$ is an ellipses, then the centers of the auxiliary circles lie on an ellipse that is polar dual to $\g$. As was mentioned  in Section \ref{sect:ell}, this follows from the Poncelet grid theorem. 
See Figure \ref{centers} for 5-periodic orbits of the outer length billiard map.
}
\end{example}

\begin{figure}[ht]
\centering
\includegraphics[width=.4\textwidth]{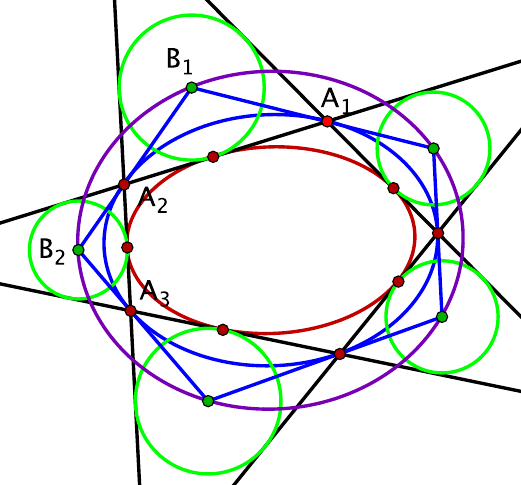}
\caption{}	
\label{centers3}
\end{figure}

Here is a more detailed explanation. Consider Figure \ref{centers3}. Point $A_2$ is the center of homothety that takes circle with the center $B_1$ to the circle with the center $B_2$. Hence the line $B_1 B_2$ bisects the angle between the common tangents $A_1 A_2$ and $A_2 A_3$. Since $A_1 A_2 A_3\ldots$ is a billiard trajectory in the middle ellipse, the line $B_1 B_2$ is tangent to this ellipse. Likewise for other points $B_i$.

Therefore the polygon $B_1 B_2\ldots$ is polar dual to the polygon $A_1 A_2 \ldots$ with respect to the middle ellipse. The latter is circumscribed about the inner ellipse, hence the former is inscribed in the ellipse, polar dual to the inner one with respect to the middle one, as claimed.

\begin{thebibliography}{99}

\bibitem{AB} A. Akopyan, A. Bobenko. {\it Incircular nets and confocal conics.} Trans. Amer. Math. Soc. {\bf 370} (2018), 2825--2854.

\bibitem{AT} P. Albers, S. Tabachnikov. {\it Introducing symplectic billiards.} Adv. Math. {\bf 333} (2018), 822--867.

\bibitem{AT1} P. Albers, S. Tabachnikov. {\it Monotone twist maps and Dowker-type theorems.} Pacific J. Math. {\bf 330} (2024),  1--24. 

\bibitem{ACT} P. Albers, A. Chavez Caliz, S. Tabachnikov. {\it Outer symplectic billiard map at infinity}. arXiv:2508.15142.

\bibitem{Av} E. Avksentiev. {\it A universal measure for a pencil of conics and the great Poncelet theorem. } Sb. Math. {\bf 205} (2014),  613--632. 

\bibitem{BBN1} L.  Baracco, O. Bernardi, A. Nardi. {\it Bialy-Mironov type rigidity for centrally symmetric symplectic billiards.} Nonlinearity {\bf 37} (2024), no. 12, Paper No. 125025, 12 pp. 

\bibitem{BB} L.  Baracco, O. Bernardi. {\it  Totally integrable symplectic billiards are ellipses.} Adv. Math. {\bf 454} (2024), Paper No. 109873, 17 pp. 

\bibitem{BBN2} L.  Baracco, O. Bernardi, A. Nardi. {\it Higher order terms of Mather's $\beta$-function for symplectic and outer billiards.} J. Math. Anal. Appl. {\bf 537} (2024), no. 2, Paper No. 128353, 20 pp. 

\bibitem{BBF} L. Baracco, O. Bernardi, C. Fierobe. {\it Starting the study of outer length billiards}.
arXiv:2504.11915.

\bibitem{BBS} S. Baranzini, M. Bialy, A. Sorrentino.
{\it Isoperimetric-type inequalities for Mather's $\beta$-function of convex billiards}.
arXiv:2509.06915

\bibitem{BM} M. Bialy, A. Mironov. {\it Integrable billiards and related topics}. arXiv:2510.03790.

\bibitem{BT} M. Bialy, S. Tabachnikov. {\it Dan Reznik's identities and more}. European J. Math. {\bf 8} (2022), 1341--1354.
 

\bibitem{DT} F. Dogru, S. Tabachnikov. Dual billiards. Math. Intelligencer {\bf 27} (2005), no. 4, 18--25.

\bibitem{EC} L. Edwards-Costa. {\it Outer length billiards on polygons}. In preparation.

\bibitem{IT} I. Izmestiev, S. Tabachnikov. {\it Ivory's theorem revisited.} J. Integrable Syst. {\bf 2} (2017), no. 1, xyx006, 36 pp. 

\bibitem{La} V. Lazutkin. {\it Existence of caustics for the billiard problem in a convex domain.} Izv. Akad. Nauk SSSR Ser. Mat. {\bf 37} (1973), 186--216. 

\bibitem{LT} M. Levi, S.  Tabachnikov. {\it The Poncelet grid and billiards in ellipses}. Amer. Math. Monthly {\bf 114} (2007), 895--908.

\bibitem{Mo} J. Moser. {\it On invariant curves of area preserving mappings of an annulus}. Nachr. Akad. Wiss. Gottingen Math-Phys. II Kl. (1962), 1--20.

\bibitem{Sa} L. Santal\'o. {\it Integral geometry and geometric probability.} Addison-Wesley Publ. Co., Reading, Mass.-London-Amsterdam, 1976. 

\bibitem{Sc} R. Schwartz. {\it The Poncelet grid.} Adv. Geom. {\bf 7} (2007), 157--175.


\bibitem{Ta95} S. Tabachnikov. {\it On the dual billiard problem.} Adv. Math. {\bf 115} (1995),
221--249.

\bibitem{Ta96} S. Tabachnikov. {\it Asymptotic dynamics of the dual billiard transformation.} J. Statist. Phys. {\bf 83} (1996),  27--37. 

\bibitem{Ta} S. Tabachnikov. {\it Geometry and billiards.} Amer. Math. Soc., Providence, RI, 2005.

\end{thebibliography}
\end{document}